\documentclass[11pt]{article}
\usepackage{fullpage}
\usepackage{makeidx}

\newenvironment{@abssec}[1]{%
     \if@twocolumn
       \section*{#1}%
     \else
       \vspace{.05in}\footnotesize
       \parindent .2in
         {\bfseries #1. }\ignorespaces
     \fi}
     {\if@twocolumn\else\par\vspace{.1in}\fi}
\newenvironment{keywords}{\begin{@abssec}{Key words}}{\end{@abssec}}
\newenvironment{AMS}{\begin{@abssec}{AMS subject classification}}{\end{@abssec}}
\newtheorem{prop}{Proposition}
\newtheorem{thm}[prop]{Theorem}
\newtheorem{lem}[prop]{Lemma}
\newtheorem{corol}[prop]{Corollary}
\newtheorem{res}[prop]{Result}
\newtheorem{dis}[prop]{Discussion}
\newtheorem{assum}[prop]{Assumptions}
\newtheorem{exmpl}[prop]{Example}
\def\bfi{\bf\textit\/}
\def\sigm{\rho}
\def\mi{{\min}}
\def\hat#1{{\widehat{#1}}}
\def\implies{\Longrightarrow}
\def\le{\leq}
\def\hatPhi{\widehat{\Phi}}
\def\hatg{\widehat{g}}
\def\hatB{\widehat{B}}

\def\tilP{{\tilde P}}
\def\tilG{{\tilde G}}
\def\hatf{\widehat{f}}

\def\calD{{\cal D}}
\def\calS{{\cal S}}
\def\hatphi{\widehat{\phi}}
\def\hatpsi{\widehat{\psi}}
\def\hatPhi{\widehat{\Phi}}

\def\ome{\omega}
\def\lam{\lambda}
\def\Lam{\Lambda}
\def\eps{\epsilon}
\def\Ome{\Omega}
\def\gam{\gamma}
\def\alp{\alpha}
\def\bet{\beta}

\def\dist{\mathop{\rm dist}\nolimits}
\def\diag{\mathop{\rm diag}\nolimits}

\def\eqbd{\mathop{{:}{=}}}

\def\trnst{{\bf T}}

\def\openC{{\rm C\kern-.48em\vrule width.06em height.6em depth-.02em 
                 \kern.48em}}
\def\openR{{{\rm I}\kern-.16em {\rm R}}}
\def\openZ{{{\rm Z}\kern-.28em{\rm Z}}}
\def\sZZ{{{\scriptstyle\rm Z}\kern-.24em{\scriptstyle\rm Z}}}
\def\openT{{{\rm T}\kern-.42em {\rm T}}}
\def\openH{{{\rm I}\kern-.16em {\rm H}}}
\def\openK{{{\rm I}\kern-.16em {\rm K}}}
\def\openL{{{\rm I}\kern-.16em {\rm L}}}
\def\openM{{{\rm I}\kern-.16em {\rm M}}}
\def\openN{{{\rm I}\kern-.16em {\rm N}}}
\def\openP{{{\rm I}\kern-.16em {\rm P}}}

\def\eqbd{\mathop{{:}{=}}}
\def\ii{{\rm i}}
\def\ee{{\rm e}}
\let\C\openC
\def\sC{{\rm C\kern-.38em\vrule width.06em height.45em depth-.02em 
                 \kern.3em}}
\let\N\openN

\let\R\openR
\let\Z\openZ
\def\eop{\hfill
        {\ \vbox{\hrule\hbox{\vrule height1.3ex\hskip0.8ex\vrule}\hrule}}
        \vskip 0.3cm \par}

\def\snc{\mathop{\rm sinc}\nolimits}

\def\spa{\mathop{\rm span}\nolimits}
\def\spec{\mathop{\rm spec}\nolimits}

\def\ran{\mathop{\rm ran}\nolimits}

\def\const{\mathop{\rm const}\nolimits}

\def\belowrightarrow#1{{{{}\over\ #1\ }\kern-1.1em\to}}
\def\l2{{L_2}}
\def\ltwo{{L_2}}
\def\norm#1{\|#1\|}
\def\ld{{L_2(\R^d)}}
\def\wsr{{W_2^s(\R^d)}}
\def\wkr{{W_2^k(\R^d)}}
\def\ws{{W_2^s}}
\def\wt{{W_2^t}}
\def\wk{{W_2^k}}
\def\Mu{{{\cal M}}}

\def\nzero{{{2\pi \Z^d}{\setminus}{0}}}
\def\Rd{\R^d}
\def\Zd{\Z^d}
\def\bks{\backslash}
\def\hatpsi{\widehat{\psi}}

\begin{document}

\title{Approximation orders of shift-invariant subspaces of $\wsr$}
\author{ Olga Holtz\thanks{Supported by the US National Science 
Foundation under Grant DMS-9872890, by Alexander von Humboldt Foundation
and by the DFG Research Center ``Mathematics for key technologies'' in
Berlin, Germany.} \\ Department of Mathematics \\
University of California \\
Berkeley, California 94720 U.S.A. \and Amos Ron\thanks{Supported by 
the US National Science Foundation under
Grants DMS-9872890, DBI-9983114. and  ANI-0085984, and by the U.S.
Army  Research Office under Contract DAAG55-98-1-0443.
} \\
Department of Computer Sciences \\
University of Wisconsin \\
Madison, Wisconsin 53706 U.S.A.}
\date{}
\maketitle



\begin{abstract} We extend the existing theory of approximation orders provided
by shift-invariant subspaces of $\l2$ to the setting of Sobolev spaces,
provide treatment of $\l2$ cases that have not been covered before, and
apply our results to determine approximation order of solutions to a 
refinement equation with a higher-dimensional solution space. 
\end{abstract}

\tableofcontents

\section{Introduction}

\subsection{General}
We are interested in this paper in the approximation order of 
{\it{shift-invariant (SI) spaces\/}} of functions defined on the Euclidean
space $\Rd$, $d\ge 1$. Such spaces play an important role in 
several areas of Real Analysis, including Spline Approximation, Wavelets,
Subdivision Algorithms, Uniform Sampling and Gabor Systems. It is not
surprising, thus, that the theory of approximation and representation from
SI spaces received significant attention and enjoyed rapid development in 
the last 10--15 years. The determination and understanding of
the {\it{approximation orders}} of these spaces  is among 
the main pillars of this theory.

As the title of this article indicates, we restrict our attention  to
approximation in Sobolev spaces: given $s\in\R$, we denote by 
$\wsr$ the 
{\bfi{Sobolev space of smoothness\/}}\index{$\wsr$ -- Sobolev space of smoothness $s$}
 $s$, i.e., the space
of all tempered distributions $f$ whose Fourier
transform\index{$\hat{f}$ -- Fourier transform of $f$} 
is locally in $\ld$ and satisfies
$$\|f\|_{\wsr}^2\eqbd \int_{\Rd}(1+|\cdot|)^{2s}|\hatf|^2<\infty.$$
(Here and elsewhere, $|\cdot|$ is the Euclidean distance in $\Rd$.)
A closed subspace $S\subset\wsr$ is {\bfi{shift-invariant}} if it is
invariant under all {\bfi{shifts}}, i.e., integer translations, or more
generally, scaled integer translations: given a fixed $h>0$,
$$\hbox{for every $\alp\in h\Zd$ and every $f\in \wsr$,}\quad
f\in S\implies f(\cdot+\alp)\in S.$$
When necessary, one identifies the underlying parameter $h$ by referring to 
$S$ as $h$-shift-invariant, and/or by denoting the SI space as $S^h$.
Also, sometimes, in order to emphasize the
ambient space $\wsr$ we write $S(\ws)$, instead of simply $S$.
The smallest SI space that contains a given $\Phi\subset \wsr$
is denoted by
$$S_\Phi\eqbd S_\Phi(\ws),$$
or, in complete detail, $S^h_\Phi(\ws)$, and we refer
then to $\Phi$ as a {\bfi{generating set}} of $S_\Phi$. The basic objective
of shift-invariant space theory is to understand the properties of
SI spaces in terms of properties of their generating sets.
In this regard we recall that an SI space generated by a {\it{singleton}}
$\Phi=\{\phi\}$ is known as {\bfi{Principal Shift-Invariant}} (PSI), while
the one generated by a {\it{finite}} $\Phi$ is referred to as
{\bfi{Finitely-generated Shift-Invariant}} (FSI).

Now, assume that we are given  a {\it{ladder}} 
$\calS\eqbd (S^h\eqbd S^h(\ws))_{h>0}$ of SI spaces. Let $k>s$. We say
that $\calS$ {\bfi{provides approximation order $k$}} (in $\wsr$), if, for
every $f\in\wkr$,
$$\dist_s(f,S^h)\eqbd \inf_{g\in S^h}\|f-g\|_{\wsr}\le C h^{k-s}\|f\|_{\wkr},$$
\index{$\dist_s$ -- distance in $\wsr$}with the constant $C$ independent of $f$ and $h$.
As is essentially known~\cite{JZ1} (and developed fully in this paper), the
above notion of approximation order depends strongly on $k$ but only
mildly on $s$.  The ladder $\calS$ is PSI or FSI if each of its
components $S^h$ is a PSI, or, respectively, FSI space.

The literature on approximation orders of SI spaces is vast, and it is not
within the scope of this paper to provide a comprehensive review of it.
We refer to the introduction and the references of \cite{BDR1} as well as to
the exposition and the references in the survey article \cite{JP}.
Many specific results on the topic are reviewed in the body of this article.
In particular, a complete characterization of the $L_2$-approximation orders
(i.e., the case $s=0$) of PSI ladders is obtained in \cite{BDR1}, while the 
analogous
results for FSI ladders are obtained in \cite{BDR1}, \cite{BDR2}
and \cite{BDR4}. There are also numerous results on approximations in other
norms, for example, in  $L_p$. Results and references in this direction
can be found in \cite{MJ1} and \cite{J}. In addition, we refer the reader 
to~\cite{HSS}, \cite{JRZ}, \cite{P}, \cite{R} for information on approximation 
properties of refinable SI spaces, and to~\cite{Meyer}, \cite{K}, \cite{J2},
 \cite{J4}, \cite{JM},~\cite{M} for results on wavelet constructions based 
on SI spaces.

\subsection{Motivation}                                                                       
While the current level of mathematical understanding of the issue of 
approximation orders of SI spaces is quite advanced, there are numerous
gaps and inconsistencies in it. This is exactly the motivation behind
the present endeavor: obtaining seamless, cohesive (and, so we hope, final)
theory. We provide a few examples for the ``gaps'' and
``inconsistencies'' in the state-of-the-art theory. Let us first
define two important classes of SI ladders: stationary ladders, and local
ones.

{\bf Definition.\/} \rm Let $\calS$ be an SI ladder. We say that
$\calS$ is {\bfi{stationary}} if, for every $h>0$,
$S^h=S^1(\cdot/h)\eqbd \{f(\cdot/h):\ f\in S^1\}$. Given a stationary ladder,
we say that $\calS$ is also {\bfi{local}} if $S^1$ is FSI and is generated
by a {\it{compactly supported}} $\Phi$.

\begin{dis} \rm

{\bf{(1)}}
Let us assume that $\calS$ is PSI, stationary and local.
Then the entire ladder is determined by the (compactly supported) generator
$\phi$ of $S^1$ (since the other spaces in the ladder are dilations of
$S^1$).
In this case, one usually refers to $\phi$ as
{\bfi{the generator of the
ladder}}. The current theory covers the case $s\ge 0$, and shows
that the approximation order in $\l2$ (as well as in $\ws$, $s>0$) provided
by such ladders is intimately related to the order of the zeros $\hatphi$, the
Fourier transform of $\phi$, has at the punctured lattice
$2\pi\Zd\bks0$ (cf.\ \S3.2). The smoothness of $\phi$, on the other hand, 
does not play any role, provided, of course, that $\phi\in \l2$ (which
is required for the definition of $\l2$-orders to make sense).
Thus, if we replace $\phi\in \l2$ by its convolution
product with a smooth generic mollifier, the $\l2$-approximation order
of the ladder, in
general, will not change.  In contrast, if $\phi\not\in\l2$ while its
Fourier transform does have the requisite zeros on $2\pi\Zd\bks0$,
the smoothing may simply result in an $\l2$-function, and the ladder
may then provide high approximation order in $\l2$ (despite the fact that 
the $\l2$-approximation order provided by the initial ladder is zero).
One expects that the extension of the notion of approximation order
to $\ws$, $s<0$ will remove the above artificial hump, and this is, indeed,
the case.

{\bf{(2)}} Retaining the same setup as in (1), it is also quite well-known
that if $\phi_1$ and $\phi_2$ are two compactly supported $L_2(\R)$-functions,
and if the PSI stationary ladder generated by $\phi_j$, $j=1,2$, provides
approximation order $k_j>0$, then the PSI stationary ladder generated by
$\phi_1\ast\phi_2$ provides (at a minimum) approximation order $k_1+k_2$.
One expects then that, if $k_2=0$, the approximation order provided
by $\phi_1\ast\phi_2$ will be at least $k_1$. This, however, is not the
case, and there are examples when the aforementioned approximation order is
{\rm{smaller}} than $k_1$. This nuisance is fixed
(in \S3.4) via the introduction of {\rm{negative}} approximation orders.

{\bf{(3)}} Let us consider now the case of local stationary FSI ladders
in $\l2$. In this case $S^1=S^1_\Phi(\l2)$, with
$\Phi\eqbd \{\phi_1,\ldots,\phi_r\}\subset\ld$ compactly supported, and with
$S^h=S^1(\cdot/h)$, $0<h<1$. A cornerstone in the analysis of the
approximation order of such ladders in the existence of a
{\rm{superfunction}}, i.e., a compactly supported function $\psi\in S^1$
whose associated local stationary PSI ladder already provides the same
approximation order as the original FSI ladder (cf.\ \S4.2). 
The existence of such a superfunction is proved in~\cite{BDR2}.
However, the Fourier transform of the
superfunction $\psi$ may vanish at the origin, a property
that denies us the existence of effective numerical approximation schemes
from its associated ladder (we refer to such superfunctions as
``bad''). In \cite{BDR4}, this problem is overcome, but
at the price of imposing an additional condition on the vector $\Phi$
(its Gramian should be invertible at the origin; cf.\ \S\ref{goodbad1} 
for a complete discussion). At the outset of the current venture, we observed
that the condition assumed in \cite{BDR4} is not {\rm{necessary}} for the
existence of a ``good'' superfunction (i.e.,
a superfunction $\psi$  for which $\hatpsi(0)\not=0$). Unfortunately, 
a good superfunction may not always exist: in \S\ref{goodbad2}
we construct an FSI vector (with $d=r=2$) for which all the compactly 
supported superfunctions are bad, dashing thereby our hope that a good 
superfunction may be proved to exist in general.

{\bf{(4)}} Let $\calS$ be a ladder as in {\bf (3)}, but assume, in addition,
that $S^1\eqbd S^1_\Phi$ is refinable, i.e., that $S^2\subset S^1$. It is
then known (see, e.g.,~\cite{Meyer},~\cite{CDM} and~\cite{R} for the PSI
case and~\cite{R} for the FSI case) that the $\l2$-approximation
orders provided by the ladder are bounded {\rm{below}} by the smoothness
of $\Phi$: if $\Phi\subset\wkr$, then the ladder provides approximation
order $k+1$ or higher. Moreover, \cite{R} proves  (for $d=1$, and
under some mild conditions on $\Phi$ for $d>1$) that approximation
order  $k+1$ is implied by the mere existence of a nonzero function 
$f$ in $S^1\cap\wkr$. However, all these results assume more than 
the smoothness of $f$ and the refinability of $S^1$: they require
in addition the entire vector $\Phi$ to lie in
$\l2$. The removal of this condition (\S4.9) leads to a conclusion
that says, essentially, that
for $\calS$ to provide some approximation order, it should contain
{\rm{one}} nonzero function of corresponding smoothness,
and nothing else.

{\bf{(5)}} Our final example still deals with refinable ladders. One
way to obtain a refinable space $S_\Phi$ is to select an $r\times r$ matrix
$P$ whose entries are trigonometric polynomials and to seek a compactly
supported vector-valued function $\Phi$ that solves the refinement
equation $\hatPhi(2\cdot)=P\hatPhi$. A major goal in this direction is
to reveal the approximation order of the stationary ladder generated
by $\Phi$ in terms of properties of $P$ 
(see~\cite{HSS}, \cite{J1}, \cite{BDR4}). 
The ultimate known result in this direction, \cite{BDR4}, 
requires a regularity condition
on $\Phi$ that necessarily fails once the above refinement equation
has (in a nontrivial way) more than one solution.  Thus, there is no
theory at present that deals with the approximation orders of refinable
vectors, once the refinement equation has multiple solutions.
\S5 deals with the approximation order of stationary refinable ladders 
and provides a novel theory for the case when multiple solutions to the 
same refinement equation exist.

\end{dis}

\subsection{Layout of this article}

In the introductory \S\ref{prelude}, we define the notions of 
shift-invariance and approximation order and make several basic 
observations that will be extensively used in the sequel.

Section~\ref{psisec} is devoted to PSI ladders. The section begins 
with a summary of the known characterization of the $\l2$-approximation 
orders provided by PSI stationary and nonstationary ladders.  These 
results are then extended to general $\ws$ (the end of \S\ref{charpsisec}) 
and connected with the notion of the Strang-Fix conditions 
(\S\ref{sfsection}) and  polynomial reproduction 
(\S\ref{polreppsi}). The results from \S\ref{sfsection} are in
turn  used in~\S\ref{conssection} to analyze the dependence of the 
approximation order notion on the value of $s$, i.e., on the space where 
the error is measured. The issue of negative approximation orders is
discussed in \S\ref{negsection}. 

FSI ladders are considered in \S\ref{fsisec}. It begins, analogously to
\S\ref{psisec}, with a summary on the $\l2$-approximation orders 
of FSI spaces and with the extension of these results to the setting of 
Sobolev spaces. This takes up \S\ref{charfsisec}. Section~\ref{supersec} 
focuses on the  notion of a superfunction, which is instrumental in the 
reduction of the FSI case to the PSI case.  This notion is further used 
in~\S\ref{polrepfsi} to understand polynomial reproduction from FSI
 spaces and in \S\ref{consFSI}  to establish the consistency of 
the notion of approximation order as we operate in different
Sobolev spaces.  However, not all superfunctions are equally useful,
as is made clear in~\S\ref{goodbad1} and~\S\ref{goodbad2}.
Regardless of whether `good' superfunctions exist in the underlying
FSI space, there is an alternative method proposed in~\S\ref{estim} 
that can always be used to estimate approximation orders. The usefulness 
of that alternative approach is demonstrated by an example in~\S\ref{c1cubic}.  
Section~\ref{smoothsec1} is devoted to refinable FSI spaces. It shows that
the approximation order of stationary refinable FSI spaces is bounded below 
by (a weak variant of) the decay rate of the Fourier transform of any
(nonzero) function in the space.

In \S\ref{refinesec}, applications of the theory from the preceding 
sections  to multiple solutions to a refinement equation are developed. 
We start by discussing, in~\S\ref{rfsec}, the structure of the 
solutions space to a refinement question. In~\S\ref{cohntsec}, we 
introduce the notion of coherent approximation orders, which bundles 
together different solutions to the same refinement equation.  
In~\S\ref{unisuper}, the notion of coherent approximation order is 
associated with a corresponding (novel) notion of universal supervectors; 
those lead to a uniform way of constructing superfunctions in all the
FSI spaces that are generated by the various solutions to the given 
refinement equation. 

For the convenience of the reader, some of the notations used in the paper 
are collected immediately before the bibliography.

\section{\label{prelude}SI ladders: the prelude}


We start our analysis with a few elementary, yet very useful,
observations concerning the interplay between approximation orders 
in $\wsr$ on the one hand, and in $\ld$ on the other hand.

As mentioned before, the symbol $\wkr$ 
denotes the {\bfi{Sobolev space of smoothness\/}} $k$. Note also
the isometry\index{$J_s$ -- canonical isometry from $\wsr$ into $\ld$} 
\begin{equation}  J_{-k}: \ld \to \wkr : 
f\mapsto \left(  (1+|\cdot|^2)^{-k/2} \widehat{f} \right)^\vee. \label{map}
\end{equation}
Recall that the Sobolev spaces are ordered
by embedding: $\ws \hookrightarrow \wt$ whenever $s\ge t$. 

\subsection{Shift-invariance}

The notion of shift-invariance is valid in any function space $F$ whose 
elements are defined on $\Rd$, and
is certainly not specific to $\wsr$. Given such a space $F$,
we consider now SI spaces that are invariant under
{\it{integer}} translations; thus, we refer to  a closed subspace
$S\subset F$ as  {\it{shift-invariant\/}} 
if $S$  is  invariant under multi-integer shifts
$$ s\in S  \Longrightarrow  s(\cdot - \alpha)\in S, \qquad \alpha \in \Z^d. $$ 
 
In agreement with the definitions of PSI and FSI ladders,
a {\bfi{principal shift-invariant (PSI)  space\/}} $S_\phi$\index{$S_\phi$ --
PSI space generated by $\phi$}
is  generated by a single function $\phi\in F$ as the closure of
$$ \spa [\phi(\cdot - \alpha) : \alpha \in \Z^d] $$
in the topology of $F$, while 
a {\bfi{finitely generated shift-invariant (FSI)  space\/}} 
$S_\Phi$\index{$S_\Phi$ -- FSI space generated by $\Phi$}
is the closure of $\sum_{\phi\in \Phi}S_\phi$, with $\Phi$
a finite subset of $F$.

It is known that an FSI subspace of $\ld$ can be characterized 
on the  Fourier domain as follows:

\begin{res}[\cite{BDR2}] \label{l2cont} For $\Phi\subset \ld$,
\begin{equation}
S_\Phi(\ld)=
\{ f\in \ld :
\widehat{f}=\tau^* \widehat\Phi, \; \tau \; {\rm measurable},\;\; 
\tau(\cdot +\alpha)=\tau, \;\;
{\rm all} \;\; \alpha \in 2\pi\Z^d \}. \label{cont1}
\end{equation}
\end{res}

That is, the Fourier transform of an element of $S_\Phi(\l2)$ is
the inner product of two vector-valued functions: the vector $\tau$
(whose entries are measurable and $2\pi$-periodic but otherwise arbitrary),
and the vector $\widehat\Phi$. Note that we tacitly assume that the entries
of $\tau$ are indexed by $\Phi$ (or by the same index set that is used
to index $\Phi$).

Since the operators $J_s$ commute with translations, one 
easily checks that
\begin{equation}
S_\Phi(\ws)=J_{-s}S_{J_s\Phi}(\l2),\label{jsiden}
\end{equation}
which, together with Result~\ref{l2cont}, leads to the following

\begin{corol} \label{sobcont}  For $\Phi\subset\wsr$,
\begin{equation} S_\Phi(\ws)=
\{ f\in \wsr : 
\widehat{f}=\tau^* \widehat\Phi, \; \tau \; {\rm measurable}, \;\; 
\tau(\cdot +\alpha)=\tau, \;\;
{\rm all} \;\; \alpha \in 2\pi\Z^d \}. \label{cont2} \end{equation}
\end{corol}


\subsection{Approximation orders}

The basic idea leading
to the notion of approximation order is very simple. It is
the heuristic understanding that increasing the density of translations
used to define an SI space may improve   their approximation ``power''.
At the same time, for numerical reasons (and also for deeper theoretical
reasons), one would, almost always, change the generator(s) of the SI
space when switching from $S^1$ to $S^h$, $h<1$: the new generators should
be more localized, and one way, sometime adequate sometime not, to
modify the generators is by dilation (see the definition of a stationary 
ladder in \S1.2).

%
%
%
%
The following, simple but important, result connects
the approximation orders of SI ladders in  $\l2$ to the analogous
approximation orders in $\ws$:

\begin{prop} \label{iniprop}
The ladder $\calS=(S^h=S^h(\ws))_h$ provides approximation 
order $k$  in $\ws$ if and only if $(J_s S^h)_h$ provides approximation 
order $k{-}s$ in $\l2$, where $J_s$ is defined as in~(\ref{map}).  \end{prop}

\sl Proof. \rm  
$J_s$ is an isometry from
$\wk$ to $W^{k-s}_2$ as well as from $\ws$ to $\l2$.
Thus, if $J_s\calS$ provides approximation $k-s$ in $\l2$ then, for
every $f\in\wk$,
$$  \dist_s(f,S^h)=\dist_0(J_sf, J_s S^h)_{\l2}\leq  C
h^{k-s} \|J_sf\|_{W_2^{k-s}}= C h^{k-s} \|f\|_{\wk}.$$ 
Hence $\calS$ provides approximation order $k$ in $\ws$. The converse
is proved in the same manner.
\eop

\medskip
As already indicated before, the
two most important cases of SI ladders are   
\begin{itemize}
\item {\bf PSI:\/} each $S^h$ is an $h$-dilate of some PSI space,
i.e., $S^h=S_{\phi_h}(\cdot/h)$; a PSI ladder may be stationary
or nonstationary depending on whether or not the generator $\phi_h$ 
of $S^h$ is independent of $h$. 

\item {\bf FSI:\/} each $S^h$ is an $h$-dilate of some FSI space
$S_{\Phi_h}$; an FSI ladder, just like a PSI ladder, may be 
stationary or nonstationary. 
\end{itemize}

Nonstationary FSI ladders are broad enough to cover almost all situations 
of interest in  applications. Thus it is of primary importance to be able 
to characterize the approximation orders provided by such ladders.
It turns out that nonstationary ladders are useful not only on their own,
but also as a tool for analyzing stationary ladders. 

\begin{corol}\label{inicorol} An FSI ladder 
$(S^h\eqbd S_{\Phi_h}(\cdot/h))_h$  provides approximation order $k$ in $\ws$ 
if and only if the  FSI ladder 
$(S_{\Psi_h}(\cdot/h))_h$, $\widehat\Psi_h\eqbd  
 (1+|\cdot/h|^2)^{s/2}\widehat\Phi_h$, provides approximation 
order $k{-}s$ in $\l2$. \end{corol}

\sl Proof. \rm In view of Proposition \ref{iniprop}, we only
need to identify $(J_s S^h)(h\cdot)$ as
$S_{\Psi_h}$,  with $\Psi_h$ as in the corollary. Now, by Corollary
\ref{sobcont},
$f\in S_{\Phi_h}(\ws)$  iff $f\in\ws$ and  $\hatf=\tau^*\hatPhi_h$, $\tau$ 
$2\pi$-periodic.  Thus, $f\in S^h$ iff $f\in \ws$ and 
$\hatf=\tau^*\hatPhi_h(h\cdot)$, with $\tau$ $2\pi/h$-periodic. Thus $f\in 
J_sS^h$ iff $f\in \l2$, and 
$$\hatf= (1+|\cdot|^2)^{s/2}\tau^*\hatPhi_h(h\cdot).$$
Dilating the last equation, we obtain that $f\in
(J_s S^h)(h\cdot)$ iff $f\in \l2$ and 
$$\hatf=(1+|\cdot/h|^2)^{s/2}\tau^*\hatPhi_h$$
for a $2\pi$-periodic $\tau$.
By Result \ref{l2cont}, this last condition is equivalent to $f$ being in 
$S_{\Psi_h}(\ltwo)$.
\eop

Note that the ladder associated with $(\Psi_h)_h$ in the above result is
nonstationary even when we assume the original one to be stationary, i.e.,
when we assume $\Phi_h$ to be independent of $h$.

\section{\label{psisec}PSI ladders}

We start our study of PSI ladders by recalling the
characterization of the $\l2$-approximation order
of these spaces. We then extend the result to the Sobolev space $\ws$.
The general result is then connected with the
notions of the Strang-Fix conditions and  polynomial reproduction. In
turn, those latter notions allow us to understand the dependence of the
approximation order notion on the value of $s$, i.e., on the space where
the error is measured.

\subsection{\label{charpsisec}Approximation orders of PSI ladders}

Note that the first part
of the next result is not entirely a special case of second part
(although it can be derived from it with ease).

\begin{res}[{\cite[Theorems~1.6 and~4.3]{BDR1}}] \label{psil2}
\hfill  
\begin{enumerate}\item The stationary PSI ladder
$\calS=(S^h\eqbd S^h(\l2))$, with $S^h=S_\phi(\cdot/h)$, $\phi \in \ld$, 
provides  approximation order $k$ if and only if there exists a neighborhood 
$\Omega$ of $0$ such that
$$ { [\widehat{\phi}, \widehat{\phi}]^0  \over [\widehat{\phi}, 
\widehat{\phi}] } \;
{1\over |\cdot|^{2k}}\in L_\infty(\Omega). $$ 
Here $[\widehat{\phi},\widehat{\phi}]\eqbd \sum_{\alpha \in 2\pi \Z^d} 
|\hat\phi(\cdot+\alpha)|^2,$\index{$[\cdot,\cdot]$ -- bracket product in $\ld$}
 $[\widehat{\phi},\widehat{\phi}]^0\eqbd \sum_{\alpha \in \nzero} 
|\hatphi(\cdot+\alpha)|^2.$\index{$[\cdot,\cdot]^0$ -- truncated bracket 
product in $\ld$}
\item The  nonstationary PSI ladder
$\calS=(S^h\eqbd S^h(\l2))$, with $S^h=S_{\phi_h}(\cdot/h)$, $\phi_h \in \ld$, 
provides approximation  
order $k$ if and only if, for some $h_0>0$ and some neighborhood $\Omega$
of $0$, the collection of functions 
$$ {[\widehat\phi_h, \widehat\phi_h]^0  \over 
 [\widehat\phi_h, \widehat\phi_h] }\; {1\over (|\cdot|^2+h^2)^k}, 
\qquad 0<h<h_0, $$
lies in $L_\infty(\Omega)$ and is bounded there.   \end{enumerate} \end{res}

Combining Proposition~\ref{iniprop} and 
Result~\ref{psil2} we obtain the analogous result for Sobolev spaces.

\begin{thm} \label{psisob} Let $s\in \R$ and $k>s$. Assume also
that $k$ is nonnegative.\hfill  
\begin{enumerate}
\item The stationary PSI ladder 
$\calS=(S^h\eqbd S^h(\ws))$, with $S^h=S_\phi(\cdot/h)$, $\phi \in \ws$,
provides approximation  order $k$ if and only if there exists a
neighborhood $\Omega$ of $0$ such that
\begin{equation}
\Mu_{\phi,s}\eqbd {  [\widehat{\phi}, \widehat{\phi}]^0_s \over
 [\widehat{\phi}, \widehat{\phi}]_s } \; {1\over |\cdot|^{2k-2s}}
\in L_\infty(\Omega). 
\label{statPSI} \end{equation} 
Here $[\widehat{\phi},\widehat{\phi}]_s\eqbd \sum_{\alpha \in 2\pi \Z^d} 
|\widehat\phi(\cdot+\alpha)|^2  
|\cdot+\alpha|^{2s}$,\index{$[\cdot,\cdot]_s$ -- bracket product
in $\wsr$} $[\widehat{\phi},\widehat{\phi}]^0_s\eqbd 
\sum_{\alpha \in \nzero} |\widehat{\phi}(\cdot+\alpha)|^2 
|\cdot+\alpha|^{2s}.$\index{$[\cdot,\cdot]^0_s$ -- truncated bracket product in $\wsr$}
\item The nonstationary PSI ladder 
$\calS=(S^h \eqbd  S^h(\ws))$, with $S^h=S_{\phi_h}(\cdot/h)$, $\phi_h \in \ws$,
provides approximation  
order $k$ if and only if, for some $h_0>0$ and some neighborhood $\Omega$
of $0$, the collection of functions 
\begin{equation} { [\hatphi_h,\hatphi_h]^0_{s}  \over 
[\hatphi_h,\hatphi_h]_{s,h}  } \;  {1\over (|\cdot|^2+h^2)^{k-s}}, \qquad 
0<h<h_0 \label{nonstatPSI} \end{equation}
lies in $L_\infty(\Omega)$ and is bounded there.
Here,
$ [\hatphi_h,\hatphi_h]_{s,h} \eqbd  [\hatphi_h,\hatphi_h]_s^0+
|\hatphi_h|^2 (|\cdot|^2+h^2)^{s}.$ 
\end{enumerate} \end{thm}

\sl Proof. \rm  The second part of the current theorem follows from 
the second part of Result~\ref{psil2} and the PSI case
of Corollary~\ref{inicorol}. Together, these two results yield the requisite
characterization, but with $[\hatphi_h,\hatphi_h]_{s}^0$ replaced by
$\sum_{\alp\in 2\pi\Zd\bks0}|\hatphi_h(\cdot+\alp)|^2(|\cdot+\alp|^2+h^2)^s$.
However, for $\alp\not=0$, we can replace 
$(|\cdot+\alp|^2+h^2)^s$ by its equivalent expression
$|\cdot+\alp|^{2s}$.

It remains to show that
in the stationary case, i.e., when $\phi_h=\phi$ for all $h$, 
(\ref{nonstatPSI}) is equivalent to (\ref{statPSI}). The fact that the former
{\it{implies}} the latter is obvious (one simply should take
$h\to 0$ in (\ref{nonstatPSI}) and invoke the uniform boundedness of the
collection of functions that appears there). For the converse, we observe
that (when $\phi_h\eqbd\phi$ for all $h$)
the uniform boundedness of the functions in (\ref{nonstatPSI})
is equivalent to the validity of the inequalities
$${[\hatphi,\hatphi]_s^0
\over |\hatphi|^2 }\le  {(|\cdot|^2+h^2)^k\over c-(|\cdot|^2+h^2)^{k-s}}
\qquad {\rm a.e.} ,$$
for some positive $c>0$.
Moreover, since we assume $k-s>0$, we can force $(|\cdot|^2+h^2)^{k-s}<c$ by 
making $h$ small enough and changing $\Ome$ if necessary. This leaves us with
$${[\hatphi,\hatphi]_s^0
\over |\hatphi|^2 }\le C (|\cdot|^2+h^2)^k \qquad {\rm a.e.}$$
as the requisite boundedness. This is definitely implied by (\ref{statPSI}),
as the left hand side in the display above is independent of $h$ and 
since $k\geq 0$.
\eop

{\bf Remark on notation.\/} For brevity, we will use 
in the sequel the expressions `$S_\phi(\ws)$ provides approximation order 
$k$' and `$\phi$ provides approximation order $k$ in $\ws$' to mean 
that the {\it stationary} ladder generated by $S_\phi(\ws)$ provides 
approximation order $k$ in $\ws$.

\subsection{Strang-Fix conditions \label{sfsection}}

Given $\phi\in \wsr$, and $k>0$,
one says that $\phi$ {\bfi{satisfies the Strang-Fix (SF)
condition of order $k$}} (\cite{SF}), if  $\hatphi$ has a zero of order
$k$ at each point $\alp\in 2\pi\Zd\bks0$.  
It is well known that the $\l2$-approximation order of a stationary
PSI ladder is closely related to  the order  of the SF condition satisfied
by the generator $\phi$ of the ladder.
To be precise, a full characterization requires a nondegeneracy condition 
on $\hatphi$ at the origin.
First, let us cite the $\l2$-result.

\begin{res}[{\cite[Theorems~1.14, 5.14]{BDR1}}] Assume that $0<\eta \leq 
|\widehat\phi|$ a.e. on some neighborhood $\Omega$ of the origin.
Let $A\eqbd  \Omega +\nzero$. If
$\widehat{\phi}\in W^\rho_2(A)$ for some $\rho>k+d/2$, then
$S_\phi(\l2)$ provides approximation 
order $k$ (in $\l2$) if and only if $\phi$ satisfies the Strang-Fix 
conditions of order $k$, i.e., near the origin
\begin{equation}
|\widehat{\phi}(\cdot +\alpha)|=O(|\cdot|^k) \qquad \hbox{\rm for all} \quad
\alpha \in \nzero. \label{SF}
\end{equation} \end{res}

Here $W^\rho_2(A)$ is the local version of $W^\rho_2(\Rd)$; cf.\
\cite[Chapter~7]{A}. For our purposes, it
is only important that the norm on $W^\rho_2(A)$ has a subadditivity
property, i.e., 
\begin{equation} \sum_\alpha \|f\|^2_{W^\rho_2(\alp+\Omega)} \leq 
\const  \|f\|^2_{W^\rho_2(A)},  \label{subadd}
\end{equation}
and that the Sobolev embedding theorem for such spaces still holds,
in particular, that  the bounded (compact) embedding
\begin{equation}
W^\rho_2(\alp+\Omega)  \hookrightarrow 
 W^k_\infty(\alp+\Omega) \label{imbed}
\end{equation} 
is valid.
Note that the condition $\widehat{\phi}\in W^\rho_2(A)$ is 
weaker than the more traditional decay condition on $\phi$
$$ |\phi|=O(|\cdot|^{-k-d-\epsilon}), \quad \epsilon>0,$$
which implies global smoothness of $\widehat{\phi}$.

We now show that the Strang-Fix conditions also characterize 
approximation power in a Sobolev space.  

\begin{thm} \label{sfthm}  Let $k\geq 0$, $s<k$, $\phi\in \ws$.
Suppose that, for some $\eta_1,\eta_2>0$ and for some ball $\Omega$ centered 
at the origin,
\begin{equation}
 \eta_1 \leq | \widehat\phi|\le \eta_2 \qquad \hbox{\rm a.e. on} \quad \Omega 
\label{techn1} 
\end{equation}
\begin{equation}
\| \widehat\phi \|^2_{k,A} \eqbd  \sum_{\beta \in \nzero} |\beta|^{2s} \max_{\gamma  :  |\gamma| \leq k}
\|D^\gamma \widehat{\phi}\|_{L_\infty(\beta+\Omega)}^2 <
\infty, \label{techn2} 
\end{equation}
where $A$ denotes the set $\Omega +\nzero$ and $D^\gamma$ denotes the
monomial derivative of order $\gamma\in \Zd_+$ normalized by 
$\gamma!$.\index{$D^\gamma$ -- normalized differential operator}
 Then $S_\phi(\ws)$ provides approximation order
$k$  (in $\ws$) if and only if~(\ref{SF}) holds. \end{thm}

\sl Proof. \rm   Set 
\begin{equation} 
R\eqbd |\cdot|^{-2s}\sum_{\beta \in \nzero} |\widehat{\phi}(\cdot +\beta)|^2 
\; |\cdot +\beta|^{2s}.
\end{equation}
Suppose $\phi$ provides approximation order $k$ in $\ws$.
Then~(\ref{statPSI}) holds by~Theorem~\ref{psisob}, or
equivalently, a.e. on $\Ome$,
$${R\over |\hatphi|^2+R}=O(|\cdot|^{2k-2s}).$$
Since $k>s$, and $\hatphi$ is bounded on $\Ome$, we conclude that, around
the origin, $R=O(|\cdot|^{2k-2s})$. This readily implies
(\ref{SF}).

Now suppose $\phi$ satisfies~(\ref{SF}).  With $R$ as above, we invoke
(\ref{techn2}) to conclude that
$$\|R\|_{L_\infty(\Omega)}\le C|\cdot|^{2k-2s}\|\hatphi\|_{k,A}^2=
O(|\cdot|^{2k-2s}).$$
However, the left-hand side $\Mu_{\phi,s}$
of~(\ref{statPSI}) equals 
$${|\cdot|^{2s-2k} R\over |\hatphi|^2+R}.$$
We have just argued that the numerator in this expression is bounded.
The denominator of the expression is bounded away of zero thanks to 
(\ref{techn1}).  This implies~(\ref{statPSI}).  
\eop

Note that condition  (\ref{techn2}) was required only for the ``if''
implication in the above result.

\begin{corol} \label{notsfcor}   In the notation of Theorem \ref{sfthm}, let
$\rho>k+d/2$ and let $A\eqbd  \Omega +\nzero$. 
Then the conclusions of Theorem \ref{sfthm} remain valid when we replace 
condition (\ref{techn2}) by: 
\begin{description}
\item{(i)} for $s\le 0$, the condition that
$\widehat{\phi}\in W^\rho_2(A)$.
\item{(ii)} for $s\ge 0$, the condition that
$(1+|\cdot|^2)^{\rho/2}\phi\in\wsr$,
or the stronger condition that $\phi\in\ws$ and
$\phi=O(|\cdot|^{-k-d-\eps})$, $\eps>0$, at $\infty$.
\end{description}
\end{corol}

Note that the first condition  in (ii) above implies, whenever
$s\ge 0$, that $\hatphi\in
W_2^\rho(\R^d)$, hence is stronger than the condition assumed in (i).

\sl Proof. \rm It is clearly sufficient to prove
that each of the conditions in (i) and
(ii) implies (\ref{techn2}). 

\rm (i): Using~(\ref{imbed}), together with the fact that that the sets
$(\bet+\Ome)$ are all translates of $\Ome$, the Sobolev embedding 
theorem applies to yield that, for $s\le 0$,
$$\|\widehat{\phi} \|_{k,A}^2 \leq  
C_1\sum_{\bet \in \nzero}  |\bet|^{2s}
\| \widehat{\phi}\|_{W^\rho_2(\Omega +\bet)}^2\le
C_2\sum_{\bet \in \nzero} 
\| \widehat{\phi}\|_{W^\rho_2(\Omega +\bet)}^2.$$
The right-hand-side in the above
is bounded, thanks to~(\ref{subadd}), by a constant
multiple of $\| \widehat{\phi}\|_{W^\rho_2(A)}^2$. Hence 
condition~(\ref{techn2}) is satisfied. 

(ii): The second condition in (ii) clearly implies the first one. Now assume
the first condition in (ii), i.e., that $f\eqbd (1+|\cdot|^2)^{\rho/2}\phi\in
\ws$. Then $\hatf$ is locally in $\l2$ and 
$$\sum_{\bet\in\nzero}|\bet|^{2s}\|\hatf\|_{L_2(\bet+2\Ome)}^2\le C
\|f\|^2_{\ws}<\infty.$$
However, $\|\hatphi\|_{W^\rho_2(\bet+\Ome)}\le
C\|\hatf\|_{L_2(\bet+2\Ome)}$, and the argument in the proof of (i) then
applies to yield (\ref{techn2}). \eop

\subsection{Approximation orders are independent of the underlying $\ws$ space  
\label{conssection} }

We are now in a position to observe that the definitions of approximation
order, if made with respect to different Sobolev spaces, are consistent
in the following sense.

\begin{prop}  \label{eqprop} If $S_\phi(\ws)$ provides approximation
order $k\geq 0$, $k>s$, then $S_\phi(\wt)$ provides the same approximation 
order for any $t\leq s$.
\end{prop}

\sl Proof. \rm  First note that $\phi$ is an element of $\wt$
whenever $t\leq s$, since $\ws$ is embedded in $\wt$.

Now note that, by Theorem~\ref{psisob},  $\phi$ provides approximation
order $k$ if and only if~(\ref{statPSI}) holds. The left-hand
side $\Mu_{\phi, s}$ of~(\ref{statPSI}) satisfies
\begin{equation}
 \left( 1- \Mu_{\phi,s}|\cdot|^{2(k-s)} \right)
\sum_{\beta \in \nzero} |\widehat{\phi}(\cdot +\beta)|^2 
|\cdot +\beta|^{2s} =\Mu_{\phi,s}|\widehat{\phi}|^2 |\cdot|^{2k}.
\label{useful} \end{equation}
Since $ \|\Mu_{\phi,s}\|_{L_\infty(\Omega)} \leq \const_{\phi,s}$,
the set $\Omega$ can be assumed to be small enough so that, e.g.,  
$$ 1- \Mu_{\phi,s}|\cdot|^{2(k-s)} \geq  1/2  \qquad
\hbox {\rm a.e. on} \quad \Omega. $$
Then 
$$ \sum_{\beta \in \nzero} |\widehat{\phi}(\cdot +\beta)|^2 
|\cdot +\beta|^{2t} \leq 
\sum_{\beta \in \nzero} |\widehat{\phi}(\cdot +\beta)|^2 
|\cdot +\beta|^{2s} \leq 2
\Mu_{\phi,s}|\widehat{\phi}|^2 |\cdot|^{2k}. $$
This implies
$$\Mu_{\phi,t}\leq 
{ \sum_{\beta \in \nzero} |\widehat{\phi}(\cdot +\beta)|^2 
|\cdot +\beta|^{2t} \over |\hatphi|^2 |\cdot|^{2k}}\leq 
2\Mu_{\phi,s}. $$
Thus, $S_\phi$ provides approximation order $k$ also in $\wt$.
\eop  
   
Proposition~\ref{eqprop} shows that $\phi$ provides approximation order on the
whole interval $\{\wt : t\leq s\}$ of Sobolev spaces once
it does so in the space $\ws$.

Let us now show that, under the regularity assumptions already used in 
Theorem~\ref{psisob}, a converse also holds. 

\begin{thm} \label{eq2prop}  Let $t<s<k$, $k\geq 0$. Suppose $\phi\in \ws$
and that it satisfies (\ref{techn1})--(\ref{techn2}) (with respect to $s$). 
Then $S_\phi(\ws)$ provides approximation order $k$  iff $S_\phi(\wt)$ provides
that same approximation order.
\end{thm}

\sl Proof.  \rm  The ``only if'' implication was proved in Proposition
\ref{eqprop} without appealing to (\ref{techn1})--(\ref{techn2}).

We prove the ``if'' assertion as follows. First, since $\phi$
provides approximation order $k$ in $W_2^t$, while satisfying
(\ref{techn1}), it must satisfy the SF conditions or order $k$
(we do not need (\ref{techn2}) for that part), by virtue of
Theorem \ref{sfthm}. Then, once $\phi$ satisfies the SF conditions
of order $k$, the facts that it belongs to $\ws$ and satisfies
(\ref{techn1})--(\ref{techn2}) imply, again by Theorem \ref{sfthm},
that it provides approximation order $k$ in $\ws$.
\eop

\subsection{Negative approximation orders \label{negsection}}

When would it make sense to have a shift-invariant space
that provides {\em negative\/} approximation order? 

Suppose we form a convolution of two compactly 
supported distributions $\phi_i$, $\widehat{\phi}_i(0)\neq 0$, $i=1,2$.  
If $\phi_i$, $i=1,2$, provides a positive approximation 
order $k_i>0$, then their convolution $\psi$
provides approximation order $k_1+k_2$. This is to be expected
if one assumes that the Strang-Fix conditions are equivalent
to approximation power (which is almost true): then
$\widehat{\psi}=\widehat{\phi_1}\widehat{\phi_2}$ 
satisfies the Strang-Fix conditions of order $k_1+k_2$ whenever
each $\widehat{\phi_i}$ satisfies the Strang-Fix conditions of order
$k_i$. A sample of a rigorous statement in this direction is as follows:

\begin{prop}  \label{conprop}  
Let $\phi_i\in W_2^{s_i}$, $i=1,2$ be nondegenerate compactly
supported distributions, i.e.,
$\widehat{\phi}_i(0)\neq 0$, and provide approximation order 
$k_i>0$, $k_i>s_i$, $i=1,2$ in their respective spaces.
Then $\psi\eqbd  \phi_1 * \phi_2$\index{$\ast$ -- convolution} 
provides approximation order $k_1+k_2$ in $W_2^{s_1+s_2}$.
\end{prop} 

\sl Proof. \rm  First observe that $\psi \in W_2^{s_1+s_2}$.
Indeed, 
\begin{equation} \|\psi\|_{W_2^{s_1+s_2}}^2 \leq 
\| |\widehat{\phi}_1| (1+|\cdot|^2)^{s_1/2} \|_\infty 
\| |\widehat{\phi}_2| (1+|\cdot|^2)^{s_2/2} \|_\infty 
\|\phi_1\|_{W_2^{s_1}}\|\phi_2\|_{W_2^{s_2}}. \label{compare}
\end{equation}
Since the $\phi_i$'s are compactly supported, their Fourier
transforms are entire functions. Moreover, the products 
$\widehat{\phi}_i (1+|\cdot|^2)^{s_i/2}$, $i=1,2$, are in $\l2$, and,
therefore,
their inverse transforms are in $\l2$, too. Those inverse
transforms are 
the result of applying a singular convolution operator to 
$\phi$. Since the convolutor decays rapidly at $\infty$, and since $\phi$
is compactly supported, the result decays rapidly at $\infty$. Altogether,
we conclude that $(\hatphi_i (1+|\cdot|^2)^{s_i/2})^\vee$ is in $L_1$, 
and consequently each $\hatphi_i (1+|\cdot|^2)^{s_i/2}$ must tend to zero at 
infinity.  Therefore their $L_\infty$-norms
must be finite. So, the right-hand side of~(\ref{compare}) is finite,
hence $\psi$ is in $W_2^{s_1+s_2}$.

Now, since we assume $\phi_i$ to provide approximation
order $k_i$, and since 
$\hatphi_i$ is bounded around the origin,  then (cf.\ the first part in the
proof of Theorem \ref{sfthm})
$$ [\widehat{\phi}_i,
\widehat{\phi}_i]_{s_i}^0  =O(|\cdot|^{2k_i}), \qquad i=1,2,$$
where we recall the notation $[g,g]_s^0
\eqbd  [g,g]_s-|g|^2 |\cdot|^{2s}$ from Theorem~\ref{psisob}.
But \begin{equation} [\widehat{\psi},\widehat{\psi}]^0_{s_1+s_2}
 \leq [\widehat{\phi}_1, \widehat{\phi}_1]^0_{s_1}
[\widehat{\phi}_2, \widehat{\phi}_2]^0_{s_2}, \label{schwarz_ineq}
\end{equation}
hence $$ [\widehat{\psi}, \widehat{\psi}]_{s_1+s_2}^0 
=O(|\cdot|^{2k_1+2k_2}).  $$
Invoking the fact that $\hatpsi(0)\not=0$, we finally conclude that 
$${[\widehat{\psi}, \widehat{\psi}]_{s_1+s_2}^0 \over 
[\widehat{\psi}, \widehat{\psi}]_{s_1+s_2}} =O(|\cdot|^{2(k_1+k_2) 
-2(s_1+s_2)}). $$ 
This, in view of Theorem~\ref{psisob}, finishes the proof. 
\eop

Now, what if, upon convolving a given distribution $\phi_1$ 
with another distribution $\phi_2$ one discovers that the approximation
order of $\phi_1 * \phi_2$ is smaller than that of $\phi_1$? Then it is 
natural to assign a negative approximation order to the 
distribution $\phi_2$. It makes little sense to define the notion
of negative approximation order in terms of the ability to approximate
functions. We choose, instead, the following technical definition, which is
consistent with the discussion so far, as well as 
with the argument used in the proof of our last result.

{\bf Definition\/.} \rm Let $s<k\le 0$, and let $\phi\in \ws$. We say that
the (stationary ladder generated by) $\phi$ provides approximation
order $k$ (in $\ws$), if, for some neighborhood $\Ome$ of the origin,
$$ \Mu_{\phi,s}\eqbd {  [\widehat{\phi}, \widehat{\phi}]^0_s \over
 [\widehat{\phi}, \widehat{\phi}]_s } \; {1\over |\cdot|^{2k-2s}}
\in L_\infty(\Omega). $$

\smallskip
Note that the above definition is consistent with the case $k>0$.
In this case, the {\it definition} of approximation order is different, but the
characterization provided in (\ref{statPSI}) of Theorem 
\ref{psisob} is exactly in the same terms.

Equipped with this last definition, we can extend  Proposition 
\ref{conprop} as follows:

\begin{prop} \label{conprop2} 
Let $\phi_i\in W_2^{s_i}(\R^d)$, $i=1,2$ provide approximation order $k_i>s_i$ 
in $W_2^{s_i}(\R^d)$, $i=1,2$. 
If the convolution product $\psi:=\phi_1* \phi_2$ lies in   
$W_2^{s_1+s_2}(\R^d)$, then it provides approximation order $k_1+k_2$ there.
\end{prop}

\sl  Proof. \rm  Since $k_i>s_i$ for $i=1,2$, it follows directly 
from the extended definition of approximation order that
$${[\hatphi_i,\hatphi_i]^0_{s_i}\over |\hatphi_i|^2 }=O(|\cdot|^{2k_i}), \qquad
i=1,2. $$ 
These two estimates, together with the inequality~(\ref{schwarz_ineq}),
imply that ${[\hatpsi,\hatpsi]^0_{s_1+s_2} \over |\hatpsi|^2}= 
 O(|\cdot|^{2k_1+2k_2})$, hence 
$$  { [\hatpsi,\hatpsi]^0_{s_1+s_2} \over |\hatpsi|^2 |\cdot|^{2s_1+2s_2} +
 [\hatpsi,\hatpsi]^0_{s_1+s_2} } \leq { [\hatpsi,\hatpsi]^0_{s_1+s_2} \over |\hatpsi|^2 |\cdot|^{2s_1+2s_2}} =O(|\cdot|^{2k_1+2k_2-2s_1-2s_2}).    $$   
This completes the proof.   \eop


\begin{corol}  \label{concorol}
Let $\phi_i\in W_2^{s_i}(\R^d)$, $i=1,2$ be compactly 
supported distributions and provide approximation order $k_i>s_i$ 
in $W_2^{s_i}(\R^d)$, $i=1,2$. Then $\phi_1* \phi_2$  provides 
approximation order $k_1+k_2$ in $W_2^{s_1+s_2}(\R^d)$.
\end{corol}

\sl Proof. \rm This fact follows from Proposition~\ref{conprop2},
since we know already from the proof of Proposition~\ref{conprop}
that the convolution of two compactly supported distributions in
$W_2^{s_i}(\R^d)$, $i=1,2$, lies in $W_2^{s_1+s_2}(\R^d)$. \eop
 

{\bf Remark.\/} In the rest of the paper, we only consider, by default,
generators $\phi$ of stationary ladders that provide approximation
order no smaller than $0$. Note that this is the case when $\phi$
is of compact support and satisfies $\hatphi(0)\not=0$.

\subsection{Polynomial reproduction \label{polreppsi}}

We restrict our attention in this subsection to local stationary PSI
ladders, and focus on the properties of the compactly supported generator
$\phi$ of the underlying ladder. To be sure, all the results here extend,
almost verbatim, to generators $\phi$ with sufficient decay at $\infty$, for
example, $|\phi|=O(|\cdot|^{-k-d-\eps})$ at $\infty$, with $k$ the investigated
approximation order and $\eps>0$.

The theory of approximation orders of local stationary
PSI ladders focuses, and rightly so, on the satisfaction of the SF conditions
(cf.\ \S\ref{sfsection}, and also the application of those conditions in 
\S\ref{conssection}).
Under the compact support assumption on $\phi$, the SF conditions are known
to be equivalent to  the polynomial reproduction property, the latter
being the subject of the current subsection. \footnote{Prior to the
publication of \cite{BR2} and \cite{BDR1}, approximation orders of
stationary PSI ladders were usually derived directly from the
polynomial reproduction property, while the SF conditions were considered to
be a technical way for the verification of polynomial reproduction. However,
as the discussion in this article clearly shows, the SF conditions
characterize
the approximation orders of the ladder even when a slow decay of the
generator $\phi$ renders the polynomial reproduction
property meaningless.} 

The connection between the SF conditions and polynomial reproduction
is classically known, and can be dated back to Schoenberg ($d=1$,
\cite{S46}), and Strang and Fix,~\cite{SF}. See also~\cite{B}.
Our approach here follows \cite{BR}. Altogether, the results of this
subsection are included for completeness, especially since the
polynomial reproduction property in the PSI case is key to 
the understanding of the more complicated polynomial reproduction property of 
FSI spaces (\S\ref{polrepfsi}), as well as the sum rules of refinable 
FSI spaces (\S\ref{unisuper}).

Suppose that $\phi$ is compactly supported. Let us first attempt to
connect the approximation orders provided  by its stationary PSI ladder to
the SF conditions. To this end, we would like to invoke Theorem
\ref{sfthm}.  This theorem requires the satisfaction of 
(\ref{techn1}) and (\ref{techn2}). 
Condition~(\ref{techn2}) is satisfied once $\phi\in\ws$, as (ii) of
Corollary \ref{notsfcor} shows. The fact that a compactly supported
distribution belongs to some $\ws$ is well-known, and follows from the
fact that it is necessarily of finite order (as a distribution). 
As to~(\ref{techn1}), since $\hatphi$ is continuous,
this condition is presently equivalent to the nondegeneracy requirement
$$\hatphi(0)\not=0.$$
Thus  we obtain the following result:

\begin{corol} \label{sfcor}
Let $\phi$ be a compactly supported distribution, and assume
that $\hatphi(0)\not=0$. Then there exists $s\in\R$ such that $\phi\in\ws$.
Moreover, the following conditions are then equivalent, for any
given $k>0$:
\begin{description}
\item{(i)} $\phi$ satisfies the SF conditions of order $k$.
\item{(ii)} The stationary PSI ladder generated by $\phi$ provides 
approximation order $k$ (in $\ws$).
\end{description}
\end{corol}

Now, recall that reproducing polynomials of total degree less
than $k$ means that
$$ \phi\ast' \Pi_{<k} \subseteq \Pi_{<k}. $$
The symbol $\ast'$\index{$\ast'$ -- semi-discrete convolution} denotes the {\bfi{semi-discrete convolution\/}}
\begin{equation} g \ast' : f\mapsto \sum_{j\in \Z^d} g(\cdot -j) f(j),
\label{semid}\end{equation}
$\Pi\eqbd \Pi(\Rd)$\index{$\Pi\eqbd \Pi(\Rd)$ -- $d$-variate 
polynomials} is the space of all $d$-variate polynomials, and
$\Pi_{<k}\eqbd  \{p\in \Pi:\ \deg p<k\}$.\index{$\Pi_{<k}$ -- 
polynomials of total degree smaller than $k$}

One way to connect the polynomial reproduction to the SF condition is via the
following variant of Poisson's summation formula, \cite{RS},  
\begin{equation} 
\phi\ast' f=\sum_{\alpha \in \Z^d} \phi\ast(e_{\alpha}f),\quad e_\alp:x\mapsto
e^{2\pi i \alp\cdot x},
\label{sdc}
\end{equation}
which is valid for every compactly supported
$\phi$ and every $C^\infty$-function $f$ (the convergence of the right-hand-side
series is in the topology of tempered distributions).
Now, for a polynomial $f$, one
easily verifies that $\phi\ast (e_{-\alp}f)=0$ iff $\hatphi$ has a zero
of order $\deg f+1$ at $\alp$.  Thus, once $\phi$ satisfies the SF conditions
of order $k$, we have that
$\phi\ast' f=\phi\ast f$ for all $f\in \Pi_{<k}.$
This establishes the sufficiency of the SF conditions, since $\phi\ast$
always maps $\Pi_{<k}$ into itself.
On the other hand, if $\phi\ast'f$ is a polynomial
of degree $<k$, then~(\ref{sdc}) shows that
\begin{equation}
\sum_{\alpha \in \Z^d \setminus\{0\}} \phi \ast (e_{\alpha} f)
\label{intsum} \end{equation}
is also a polynomial of degree $<k$. This is possible~\cite[Proof of 
(2.10)~Lemma]{BR} only if all the summands in~(\ref{intsum}) vanish.
In conclusion, $\phi$ satisfies the SF conditions of order $k$ if and only if
$$\phi\ast'=\phi\ast,\quad \hbox{on $\Pi_{<k}$.}$$
Since, as we already said,
$\Pi_{<k}$ is an invariant subspace of $\phi\ast$ (with or
without the SF conditions), we finally need only to guarantee that
$\phi\ast$ be injective on polynomials, or equivalently, we need
to assume that $\hatphi(0)\not=0$.
Indeed, the condition $\hatphi(0)\not=0$ is necessary and sufficient
for $\phi\ast$ to be an automorphism on $\Pi_{<k}$ (for any positive integer
$k$), and we arrived at:

\begin{thm} \label{repthm}
Let $\phi$ be any compactly supported distribution with
$\widehat{\phi}(0)\neq 0$, and let $k$ be a positive integer.
Then $\phi$ provides approximation order 
$k$ in some Sobolev space $\ws$, $s<k$, if and only if
it reproduces polynomials
of total degree less than $k$. \end{thm}

\sl Remark. \rm 
As alluded to before in a footnote,
Theorem~\ref{repthm} could also be proved directly, avoiding the
use of the Strang-Fix conditions and constructing instead a quasiinterpolant 
$Q : \ws  \to S(\phi)$ such that $Sp=p$ for any $p\in \Pi_{<k}$; for a 
detailed  discussion of this method see~\cite[Section 4]{BR}.

\section{FSI ladders \label{fsisec}}

We start this section, just like in the PSI case, by recalling the
characterization of the $\l2$-approximation orders of FSI spaces
and extending the result to the setting of Sobolev spaces. We then 
focus in \S\ref{supersec}  on the notion of a superfunction, 
which  leads to the reduction of the FSI case to the PSI case. 
Besides, this notion proves to be very helpful in understanding 
polynomial reproduction from FSI spaces (see \S\ref{polrepfsi}). 
In our setting of $\ws$, it also helps to establish the consistency of 
the notion of approximation order as we vary $s$ (see \S\ref{consFSI}). 
However, not every superfunction can be used for these and/or for
other purposes, and this brings one to the
notions of `good' and `bad' superfunctions that are discussed 
in \S\ref{goodbad1}. We show, in \S\ref{goodbad2},
that there exist FSI spaces that do not contain any good superfunctions.
Regardless of whether or not good superfunctions are around, an alternative 
approach,  that is presented in \S\ref{estim}, can always be used to bound 
the approximation order from {\it above}.
The efficacy of this method is demonstrated in
\S\ref{c1cubic}, where we recover the well-known example of
$C^1$-cubics on a three-directional mesh, \cite{BH}: this is an
example of a bivariate stationary local FSI ladder  that, while
reproducing all polynomials in $\Pi_3$, fails to provide the ``expected''
approximation order $4$.  Finally, \S\ref{smoothsec1} applies the results
obtained in this section to the case when the vector $\Phi$ is 
refinable: in establishes a lower bound on the approximation order 
provided by $\Phi$ in terms of the decay of the Fourier 
transform of any non-zero function in $S_\Phi$.

\subsection{Characterization of approximation power \label{charfsisec}}

The first three results of this section form a summary of the 
known characterization of approximation power valid in $\l2$,
while the rest constitutes the characterization in the more
general setting of $\ws$. 

\begin{res}[{\cite[Theorem~2.2]{BDR4}}] \label{fsithm0}
The stationary FSI ladder
$\calS=(S^h\eqbd S^h(\l2))$, with $S^h=S_\Phi(\cdot/h)$, $\Phi \subset \ld$, 
 provides approximation order $k$ if and only if 
there exists a neighborhood $\Omega$ of $0$ such that
$$ \left( 1- \widehat{\Phi}^*G_\Phi^{-1} \widehat{\Phi} \right)
{1\over |\cdot|^{2k}}\in L_\infty(\Omega). $$ 
Here \begin{equation} 
G_\Phi\eqbd \sum_{\alpha \in 2\pi \Z^d}\widehat{\Phi}(\cdot+\alpha) 
{\widehat{\Phi}}^*(\cdot+\alpha)=\left([\widehat{\phi}, \widehat{\varphi}]
\right)_{\phi, \varphi\in \Phi}.\label{gramian}  \end{equation}
\index{$G_\Phi$ -- Gramian in $\ld$}%
Also, the expression $G_\Phi^{-1} \widehat{\Phi}$ is taken to mean
any  solution to the equation $G_\Phi \tau=\widehat{\Phi}$.
A simple linear-algebraic argument shows that the latter equation 
is always solvable whether or not $G_\Phi$ is invertible, since one 
of the rank-one terms in~(\ref{gramian}) is $\widehat{\Phi}
\widehat{\Phi}^*$.  \end{res}

\begin{res}[{\cite[Theorem~2.7]{BDR4}}] 
\label{fsithm1}
An FSI nonstationary ladder
$\calS=(S^h\eqbd S^h(\l2))$, with $S^h=S_{\Phi_h}(\cdot/h)$, $\Phi_h \subset \ld$, 
provides approximation order $k$ if and only if,
for some $h_0>0$ and some neighborhood $\Omega$ of the origin, the 
collection of functions 
$$ \left( 1- \widehat{\Phi_h}^*G_{\Phi_h}^{-1} \widehat{\Phi_h} \right)
{1\over (|\cdot|^2 +h^2)^k}, \qquad h<h_0$$
lies in $L_\infty(\Omega)$ and is bounded there. \end{res}

We also require the following equivalent formulation of the last
characterization, in which we use the notation:
\begin{equation} 
G_\Phi^0\eqbd \sum_{\alpha \in 2\pi \Z^d\bks0}\widehat{\Phi}(\cdot+\alpha) 
{\widehat{\Phi}}^*(\cdot+\alpha)=G_\Phi-
\widehat{\Phi} {\widehat{\Phi}}^* .\label{gramian0}  \end{equation}
\index{$G_\Phi^0$ -- truncated Gramian in $\ld$}

\begin{res}(another version of Result~\ref{fsithm1}) 
\label{fsithm2} The FSI nonstationary ladder $\calS=(S^h\eqbd S^h(\l2))$, 
with $S^h=S_{\Phi_h}(\cdot/h)$, $\Phi_h \subset \ld$,
provides  approximation order $k$ if and only if 
the collection of functions $(\Mu_{\Phi_h, s,h} : 0<h<h_0)$, where
$$  \Mu_{\Phi,s,h} : \omega \mapsto 
{1 \over (|\omega|^2+h^2)^k }  \inf_{v\in\sC^{\Phi}} {v^* G^0_{\Phi}(\omega) v \over 
v^* G_{\Phi}(\omega) v }   $$
is bounded in  $L_\infty(\Omega)$ for some neighborhood $\Omega$ of the
origin and some $h_0>0$.  
\end{res}

Using these results, one obtains the following characterization
of approximation power in $\ws$.

\begin{thm} \hfill  \label{fsisob}
\begin{enumerate} \item
An FSI stationary ladder $\calS=(S^h\eqbd S^h(\ws))$, 
with $S^h=S_{\Phi}(\cdot/h)$, $\Phi \subset \ws$,
provides approximation order $k> 0$ if and only if 
there exists a neighborhood $\Omega$ of $0$ such that the function
\begin{equation} \Mu_{\Phi,s} : \omega \mapsto {1\over |\omega|^{2k-2s}} 
\inf_{v\in \sC^\Phi} 
{v^* G^0_{\Phi,s}(\omega) v \over v^* G_{\Phi,s}(\omega) v } \quad
\hbox{\rm belongs to} \quad L_\infty(\Omega). \label{lstat} \end{equation}  
Here  \begin{eqnarray*} 
G_{\Phi,s} & \eqbd  & \sum_{\alpha \in 2\pi \Z^d} 
\widehat{\Phi}(\cdot +\alpha) \widehat{\Phi}^*(\cdot +\alpha) |\cdot +\alpha|^{2s}, \\
 G^0_{\Phi,s} & \eqbd  & \sum_{\alpha \in \nzero}
\widehat{\Phi}(\cdot +\alpha) \widehat{\Phi}^*(\cdot +\alpha) 
|\cdot +\alpha|^{2s}. \end{eqnarray*}\index{$G_{\Phi,s}$ -- Gramian
in $\wsr$}\index{$G_{\Phi,s}^0$ -- truncated Gramian in $\wsr$}
\item  An FSI nonstationary ladder $\calS=(S^h\eqbd S^h(\ws))$, with 
$S^h=S_{\Phi_h}(\cdot/h)$, $\Phi_h \subset \ws$, 
provides  approximation order $k\geq 0$ if and only if 
there exists a neighborhood $\Omega$ of $0$ and $h_0>0$
such that the collection of functions $(\Mu_{\Phi,s,h}: 0<h<h_0)$, with
\begin{equation}
 \Mu_{\Phi_h,s,h} : \omega \mapsto  {1\over (|\omega|^2+h^2)^{k-s} }
 \inf_v {v^*  G_{\Phi_h,s,h}^0(\ome) v \over v^*  G_{\Phi_h,s,h}(\ome) v},
\qquad \hbox{\rm is bounded in} \quad  L_\infty(\Omega). \label{lnonstat}
\end{equation}
Here,
\begin{eqnarray*}
G_{\Phi_h,s,h}(\ome) & \eqbd  & \sum_{\alpha \in 2\pi \Z^d}
\widehat{\Phi_h}(\omega +\alpha)\widehat{\Phi_h}^*(\omega +\alpha)
(|\omega +\alpha|^2 +h^2)^s,  \\
G_{\Phi_h,s,h}^0(\ome) & \eqbd  & \sum_{\alpha \in \nzero}
\widehat{\Phi_h}(\omega +\alpha)\widehat{\Phi_h}^*(\omega +\alpha)
(|\omega +\alpha|^2 +h^2)^s.
\end{eqnarray*} \end{enumerate} \end{thm}

\sl Proof.  \rm The proof is analogous to that of Theorem~\ref{psisob}.
In particular, part~2 of the current theorem is a direct consequence of 
Result~\ref{fsithm2} and Proposition~\ref{iniprop}.

Now we use the result of part 2 to derive part 1.
In the stationary case, $\Phi_h=\Phi$ for all $h$, so the left-hand side
of~(\ref{lnonstat}) becomes
$$ \Mu_{\Phi,s,h}(\omega)={1\over (|\omega|^2+h^2)^{k-s} }
 \inf_v {v^*  \sum_{\alpha \in  \nzero}
\widehat{\Phi}(\omega +\alpha)\widehat{\Phi}^*(\omega +\alpha)
(|\omega +\alpha|^2 +h^2)^s\,  v \over 
v^*  \sum_{\alpha \in 2\pi \Z^d} \widehat{\Phi}(\omega +\alpha)
 \widehat{\Phi}^*(\omega +\alpha) 
(|\omega +\alpha|^2+h^2)^s\,  v }. $$
Since the numerator of the infimum expression is bounded above and below by 
positive multiples of $v^* G^0_{\Phi,s} v$, 
the collection $(\Mu_{\Phi,s,h})$ is bounded in 
$L_\infty(\Omega)$ if and only if
the collection  of functions 
\begin{equation} \omega\mapsto{1\over (|\omega|^2+h^2)^{k-s} }
 \inf_v {v^* G^0_{\Phi,s}(\omega) v \over 
  v^* G^0_{\Phi,s}(\omega) v + v^*  \widehat{\Phi}(\omega)
 \widehat{\Phi}^*(\omega)(|\omega|^2+h^2)^s  v } \label{auxil} \end{equation}
is bounded in $L_\infty(\Omega)$. 
Since $k$ and $k-s$ are nonnegative, for a fixed $v$ and a fixed 
$\omega\in \Omega$ (assuming $\Omega$ is sufficiently small), 
the expression 
$$ {1\over (|\ome|^2+h^2)^{k-s} }
{v^* G^0_{\Phi,s}(\ome) v \over 
  v^* G^0_{\Phi,s}(\ome) v + v^*  (\widehat{\Phi}
 \widehat{\Phi}^*)(\omega) (|\ome|^2+h^2)^s  v } $$
monotonically increases  (as $h\to 0$)
to $$  {1\over |\omega|^{2(k-s)}}  {v^* G^0_{\Phi,s}(\omega) 
v \over v^* G_{\Phi,s}(\omega) v},  $$
hence~(\ref{auxil}) monotonically increases to the function in (\ref{lstat}). 
Therefore,  the collection
$(\Mu_{\phi,s,h})$ is bounded if and only if~(\ref{lstat}) holds.
\eop

{\bf Remark on notation.\/} As in the PSI case, we shall 
use in the sequel language such as `an FSI space $S_\Phi(\ws)$ provides 
approximation order $k$' and even `$\Phi\subset \ws$ provides approximation 
order $k$' and will mean by that that the FSI stationary ladder generated by 
$S_\Phi(\ws)$ provides approximation order $k$ in $\ws$.  

The $\l2$-characterizations above (Results~\ref{fsithm0}--\ref{fsithm2}) 
are connected to the notion of superfunctions. We will now discuss this 
notion, and extend it to the setting of Sobolev spaces. 

\subsection{\label{supersec}Superfunctions}

Let $S$ be an SI space, and let $\calS=(S^h\eqbd S(\cdot/h))$
be the associated stationary ladder.
A function $g\in S\subset \l2$ is a {\bfi{superfunction\/}} in $S$ 
if the PSI stationary ladder it generates provides the same approximation 
order as that of $S$ (or, more precisely, of $\calS$). For sure, $\l2$ in 
this definition can be replaced by any Sobolev space $\ws$.

The question of existence of superfunctions in FSI spaces can be answered
in the affirmative using Theorem~\ref{fsisob}.

\begin{thm} \label{super}
Any FSI space $S_\Phi \subset \ws(\R^d)$ contains a superfunction. 
\end{thm}

\sl Proof. \rm Let $\Ome$ be as in Theorem \ref{fsisob}.
 For each fixed $\omega\in \Omega$, there exists a vector
$v_0(\omega)\in\C^\Phi$ of (e.g., Euclidean) norm $1$ that minimizes the ratio 
$v^* G_{\Phi,s}^0(\omega) v/ v^* G_{\Phi,s}(\omega) v$.
Now extend $v_0$ to the cube $[-\pi, \pi)^d$ in an arbitrary way,
provided the norm $v_0$ is everywhere equal to $1$. Finally, extend
$v_0$ so defined to a $2\pi$-periodic vector-valued function.

Now, suppose that $\Phi$ provides approximation order $k$. Then,
in view of Theorem  \ref{fsisob}, the vector $v_0$ satisfies
$$ v_0 G_{\Phi,s}^0 v_0 \leq \const |\cdot|^{2k-2s} ( v_0 G_{\Phi,s}^0 v_0
+v_0^* \widehat\Phi \widehat\Phi^* v_0  |\cdot|^{2s} ) \qquad
 \hbox{\rm a.e. in } \Omega  $$ 
or, equivalently 
$$ (1-\const |\cdot|^{2(k-s)}) v_0 G_{\Phi,s}^0 v_0 \leq \const 
|\cdot|^{2k} v_0^* \widehat\Phi \widehat\Phi^* v_0 \qquad
 \hbox{\rm a.e. in } \Omega. $$
By changing $\Ome$ if needed (and using the fact that $k>s$), we obtain that
$$  v_0 G_{\Phi,s}^0 v_0 \leq C 
|\cdot|^{2k} v_0^* \widehat\Phi \widehat\Phi^* v_0 \qquad 
\hbox{\rm a.e. in } \Omega   $$
for some constant $C$. This implies that, for almost every fixed
$\ome\in\Ome$, the smallest eigenvalue of
the measurable Hermitian matrix
$H(\ome)\eqbd   (G_{\Phi,s}^0-C|\cdot|^{2k}\widehat\Phi \widehat\Phi^*)(\ome)$ 
is nonpositive. 
By Lemma~2.3.5 from \cite{RSh}, we can define a map $w$ on $\Ome$ such
that (i) for almost every $\ome\in \Ome$, $w(\ome)$ is a normalized 
eigenvector of $H(\ome)$ that corresponds to the minimal eigenvalue, and 
(ii) $w$ is measurable on $\Ome$. Without loss, we assume that our original
$v_0$ coincides with $w$ on $\Ome$. In particular, $v_0$ is now known
to be measurable.

Let $\phi$ be the (scalar) distribution whose Fourier transform satisfies
$\widehat{\phi}=v_0^*\widehat{\Phi}$. To show that $\phi$ is a superfunction
for $S_\Phi$, we only need to verify that it belongs to $S_\Phi$,
since it follows directly from the construction of $\phi$ that the space 
$S_\phi$ provides approximation order $k$. Since $\hatphi=v_0^\ast\hatPhi$,
we only need, in view of Corollary~\ref{sobcont},
to show that $\phi\in \ws$.
This final result is a simple consequence of the representation
$\hatphi=v_0^\ast\hatPhi$, using the facts that $\Phi\subset \ws$ and
that the entries of $v_0$ are bounded.
\eop

This theorem extends the known $\l2$ result. However, in the $\l2$-case, 
a superfunction was originally constructed as the orthogonal projection
$P_\Phi:\l2\to S_\Phi$ of the sinc-function
 $$ \snc(x) \eqbd  \prod_{i=1}^d {\sin(\pi x(i))
\over \pi x(i) }.$$
The fact that $P_\Phi (\snc)$ is a superfunction
follows from the general principle: 

\begin{res}[{\cite{BDR1}}] \label{supthm}
Let $S_\Phi(\l2)$ be an FSI space that provides approximation order $k_1 
\geq 0$, and let $S_g(\l2)$ be a PSI space that provides approximation order 
$k_2 \geq 0$. Then the PSI space generated
by the orthogonal projection $P_\Phi g$  on $S_\Phi$ in $\l2$ provides 
approximation order $ \min \{ k_1, k_2\}$.
\end{res}

If, in 
Result~\ref{supthm} above, we choose $g$ such that $k_2\ge k_1$, and
$S_\Phi(\l2)$ does not provide an approximation order greater than $k_1$,
then $P_\Phi g$ is a superfunction.
It is easily checked that the approximation order provided by
the space $S_{\snc}$ is infinite (i.e., exceeds any finite
$k$), hence the following corollary.

\begin{res}[{\cite{BDR1}}] \label{sincthm} An FSI space $S_\Phi(\l2)$ 
provides the same  approximation order as the PSI space generated by 
$P_\Phi (\snc)$ in $\l2$.
\end{res}

To find the $\l2$-projection on a FSI space of a function $f$, one
solves equation~(\ref{project}).

\begin{res}[{\cite{BDR1}}] \label{conthm} The $L_2$-projection 
$P_\Phi (f)$ of $f\in \l2$ on an FSI space $S_\Phi$ satisfies
\begin{equation}
(P_\Phi (f))^\wedge=\tau^*_f \widehat{\Phi} 
\label{project}
\end{equation} 
with $\tau_f$ {\it any} solution of 
$$G_\Phi\tau_f=[\widehat{\Phi},\widehat{f}].$$ 
Here, $\tau_f$ is a vector-valued function (indexed by $\Phi$) whose 
entries are measurable and $2\pi$-periodic, and the symbol $[\hatPhi,
\hatf]$ stands for $([\hatphi,\hatf])_{\phi\in\Phi}$, where
$[\widehat{f},\widehat{g}]\eqbd \sum_{\alpha\in \Z^d} 
\widehat{f}(\cdot+2\pi\alpha) \overline{\widehat{g}} (\cdot+2\pi \alpha)$.
\index{$[\cdot,\cdot]$ -- bracket product in $\ld$}

\end{res}

\noindent {\bf Remark:\/} Results~\ref{supthm},~\ref{sincthm} 
and~\ref{conthm} are all corollaries to Theorem~3.3 of~\cite{BDR1}. 

The above results extend easily to Sobolev spaces.
Indeed, Results~\ref{supthm}
and~\ref{sincthm} require only one assumption, viz.
\begin{equation}
P_{A} P_g=P_{P_A g} P_g,\label{commute}
\end{equation}
where $A$ denotes an arbitrary SI subspace of $\ws$, 
and  $P_A,P_g,P_{P_A g}$ are the  orthogonal projectors from $\ws$ onto 
$A, S_g(\ws),S_{P_A g}(\ws)$, respectively.

Under this condition, the analysis from~\cite[Section~3]{BDR1} leading
to Results~\ref{supthm} and~\ref{sincthm} goes 
through verbatim.  Since the  $\ws$-version of
(\ref{commute}) is a simple consequence of the $\l2$-version when combined
with the identity
 $P_{J_t A} J_t=J_t P_A$ (where $A$ an SI subspace of $\ws$ and $P_{J_t A}$ 
the orthogonal projector onto the space $J_t A$ in $\l2$),  we obtain
the following extension.

\begin{thm}\label{3.3}
Let $S_\Phi(\ws)$ be an FSI space that provides approximation order 
$k_1\geq 0$, and let $S_g(\ws)$ be a PSI space that provides approximation 
order $k_2\geq 0$.  Set $\psi\eqbd P_\Phi g$,
with $P_\Phi$ the orthogonal projection  of
$\ws$ onto $S_\Phi(\ws)$. Then the stationary PSI ladder
generated by $\psi$ provides approximation order $\min \{ 
k_1, k_2\}$. Specifically, for every $f\in \ws$ and $h>0$,
$$\dist_s(f,S_{\psi}^h(\ws))\le
\dist_s(f,S_\Phi^h(\ws))+2\dist_s(f,S_g^h(\ws)).$$
In particular, $S_\Phi(\ws)$ provides the same approximation 
order as the PSI space generated by $P_\Phi (\snc)$ in $\ws$.
\end{thm}

Superfunctions are obviously useful if one wishes to approximate
functions from a given SI space S, for if a superfunction $\phi$
is known explicitly, one can instead approximate from the simpler
space $S_\phi$. In addition, it is already well established in the 
$\l2$-theory that superfunctions give rise to quasi-interpolants, 
i.e., bounded linear maps into the underlying SI space $S$ that 
reproduce polynomials contained in $S$ (see, e.g.,~\cite{B},~\cite{BR}). 
Superfunctions were also used in~\cite{BDR1},~\cite{BDR4} as purely 
theoretical tools.

The natural expectation is that superfunctions play similar roles 
in the setting of Sobolev spaces. That turns out to be the case. 
In particular, the superfunction method allows us to lift painlessly various 
results from the PSI setup to the FSI one. This includes the discussion
concerning  the consistency of the definitions of approximation
orders in different Sobolev spaces, which we embark on in the next 
subsection. 


\subsection{\label{consFSI}Approximation orders are independent of 
the underlying $\ws$ space} 

\begin{prop}  \label{eqpropFSI} If an FSI space $S_\Phi(\ws)$ provides 
approximation order $k\geq 0$, $k>s$, then $S_\Phi(\wt)$ does so for any 
$t\leq s$.
\end{prop}

\sl Proof. \rm  Assuming that $S_\Phi(\ws)$ provides approximation
order $k$, Theorem~\ref{super} ascertains that $S_\Phi(\ws)$ contains a PSI
subspace $S_\phi(\ws)$ that already provides approximation order $k$.
Since $S_\Phi(\wt)\supset S_\Phi(\ws)$ (as easily follows from
Corollary~\ref{sobcont}), $\phi\in S_\Phi(\wt)$, too.
By Proposition~\ref{eqprop}, $S_\phi(\wt)$ then provides approximation
order $k$ in $\wt$, therefore $S_\Phi(\wt)$ provides approximation order
(at least) $k$ in $\wt$. \eop

Similarly to the PSI case, a converse
also holds under some regularity assumptions on the superfunction.


\begin{prop} \label{eq2propFSI}
  Let $k>s>t$, $k\geq 0$. Suppose $\phi,\Phi \subset \ws$. Suppose that
  $S_\Phi(\wt)$ provides approximation order $k$ in $\wt$, and that 
  $\phi\in S_\Phi(\wt)$ is a corresponding superfunction.  If $\phi$
  satisfies~ (\ref{techn1})--(\ref{techn2}),
then $S_\Phi(\ws)$ provides approximation order $k$ in $\ws$, too.
\end{prop}

\sl Proof.  \rm 
We apply Theorem~\ref{eq2prop}
to the function $\phi$ to show that $S_\phi(\ws)$ provides
approximation order $k$, which implies that $S_\Phi(\ws)$
also provides approximation order (at least) $k$. \eop  

The theorem highlights a central point: it is useful to know that 
an FSI space contains a ``good'' superfunction. In the current context
``good'' in interpreted as ``satisfying (\ref{techn1})--(\ref{techn2})''.
We will come back to this issue later, but first
we show how the superfunction method reduces the polynomial
reproduction issue in the FSI setup back to the simpler PSI setup.

\subsection{\label{polrepfsi}Polynomial reproduction}

Let $\Phi$ be a vector of compactly supported elements in
$\ws$, $s\in\R$. Suppose that $\Phi$ provides approximation order 
$k> 0$ in $\ws$. Let us assume, further, that $S_\Phi(\ws)$ contains
a {\it good} superfunction $\psi$ is the sense that:
\begin{enumerate}
\item $\hatpsi(0)\not=0$, and
\item $\psi$ is a {\it finite} linear combination of the shifts of
$\Phi$ (hence, in particular, is compactly supported).
\end{enumerate}
We note that the current notion of ``good'' is stronger (i.e.,
implies) the one that was discussed at the end of the last subsection
(as the argument in \S\ref{polreppsi} shows).

By our assumptions here,
$\hatpsi=v^\ast\hatPhi$, with $v$ a vector of {\it trigonometric
polynomials}. Therefore, with $(a_\phi)_{\phi\in\Phi}$ the Fourier
coefficients of the entries of $v$, we have the representation
$$\psi=\sum_{\phi\in\Phi}\phi\ast'a_\phi,$$
and each $a_\phi:\Zd\to\C$ is finitely supported. Here, $\ast'$
is the semi-discrete convolution, (\ref{semid}).

Next, since $\hatpsi(0)\not=0$, and $S_\psi(\ws)$ provides approximation
order $k$, we conclude from  Theorem~\ref{repthm} that
$\psi\ast'$ maps $\Pi_{<k}$ onto itself. Writing $\psi\ast'$ in terms
of $\Phi$ we obtains
$$\psi\ast'f=(\sum_{\phi\in\Phi}\phi\ast'a_\phi)\ast'f=
\sum_{\phi\in\Phi}\phi\ast'(a_\phi\ast' f).$$
The above representation leads to several conclusions that we summarize
in our next result:

\begin{corol} \label{fsipol} Let $\Phi$ be a compactly supported vector that 
provides approximation order $k>0$ in $\ws$, and assume that 
$S_\Phi(\ws)$ contains a good superfunction in the above sense. Then there 
exist finitely supported sequences $a_\phi:\Zd\to\C$, $\phi\in\Phi$ such 
that, for every $f\in\Pi_{<k}$, 
$$Tf\eqbd \sum_{\phi\in\Phi}\phi\ast'(a_\phi\ast f)$$
is a polynomial.  Here, $a_\phi\ast f$ denotes the {\sl discrete\/} 
convolution of
$a_\phi$ and $f_{|_{\Zd}}$.\hfill\\ \noindent
The polynomial $Tf$ is identical to the result of the
following {\it 
continuous} convolution
\begin{equation}
\sum_{\phi\in\Phi}\phi\ast (a_\phi\ast f)=\sum_{\phi\in\Phi}
\sum_{j \in\Zd}a_\phi(j )(\phi\ast f)(\cdot-j ).\label{conv}
\end{equation}
Moreover, the map $T_k\eqbd T_{|_{\Pi_{<k}}}$ is an automorphism.
\end{corol}

There are several immediate conclusions that can be derived directly
from the above corollary. For example, since  $T_k$ can be extended to a
convolution operator, it commutes with differentiation in the sense that 
$D^\gam T=TD^\gam$ for every $\gam\in\Zd_+$, and commutes also
with translations. 

A simpler consequence is as follows: since $T_k$ is an automorphism,
then every monomial $()^\alp$, $|\alp|<k$, lies in its range.
(Here, the symbol $()^\alp$ stands for the normalized monomial
$$()^\alp: x\mapsto (x)^\alp\eqbd  x^\alp/\alp!.$$\index{$()^\alpha$ -- 
normalized monomial}
We also use in the sequel $D^\alp$ for the {\it normalized} monomial
derivative.) Thus, the following is true:

\begin{corol}\label{fsipol1} Let $\Phi$ be as in Corollary \ref{fsipol}. 
Then there exist polynomials $(g_\alp)_{\alp\in\Zd_+}$ 
such that, for $|\alp|<k$,
$$\sum_{\phi\in\Phi}\phi\ast'(a_\phi\ast g_\alp)=()^\alp.$$
\end{corol}

The result shows that every $()^\alp$ is writable as
$\sum_{\phi}\phi\ast' f_{\phi,\alp}$, for suitable polynomials
$(f_{\phi,\alp})_\phi$. However, the result shows more: it decomposes
each $f_{\phi,\alp}$ into $a_\phi\ast g_\alp$, with the first factor
independent of $\alp$ (and is finitely supported), and the second independent
of $\phi$ (and is a polynomial).

The reader might wonder how realistic the assumption about the existence
of  a good superfunction is. We discuss that issue in this section as well as
in \S5. A sufficient condition for the existence of a good superfunction
as above is the invertibility, in a suitable sense,
of the Gramian $G_{\Phi,s}$ around the origin.
We also note that our results here recover
the results of~\cite[Section~3]{CHM} (cf. also~\cite{CHM1}, \cite{CHM2}).
The underlying assumption in~\cite{CHM} is that the shifts of the distributions
 $\phi\in\Phi$ are linearly independent, a condition that is significantly 
stronger than the Gramian  invertibility that we have alluded to above.
At the same time, our derivation here is simpler due to the superfunction
approach. 

Next, one might also wonder how to invert the operator $T$, i.e.,
how to compute the above polynomials $(g_\alp)_\alp$.
That inversion is the key for the so-called quasi-interpolation approach, 
and is discussed in detail in \cite{B} and \cite{BR} (in the PSI context; our 
superfunction approach already reduced the problem to that setup). At base, 
we seek a simple linear functional $\mu$ such that $\mu\ast$ inverts on 
$\Pi_{<k}$ either the convolution $\psi\ast$ or the map
$f\mapsto f\ast'\psi$. 

Among the various methods, we describe a general 
recursive approach (see \cite{CJW}, \cite{B}, \cite{CHM}).
To this end, we need first to present this approach in the 
nondegenerate PSI case, i.e., when the (single) generator 
$\psi$ satisfies the condition $\widehat{\psi}(0)\neq 0$. The superfunction
method will allow us then to lift the result to the FSI setup.

\begin{prop} \label{polpsi} Let $\psi$ be a compactly supported distribution
with $\widehat{\psi}(0)=1$ that provides approximation order $k$
in some $\ws$, $s\in\R$.  Define the polynomials $g_\alpha$,
$\alp\in\Zd_+$, $|\alp|<k$, by the recurrence 
\begin{equation} g_\alpha\eqbd  ()^\alpha -
\sum_{\beta<\alpha} c(\alp-\bet)\,g_\beta,
\label{polrec}
\end{equation}
where
$$c(\gamma)\eqbd (\psi\ast' ()^\gamma)(0)=(\psi\ast ()^\gamma)(0)=
(()^{\gamma}\ast'\psi)(0),\quad \gamma\in\Zd_+,\ |\gamma|<k.$$
Then these polynomials satisfy 
\begin{equation}
()^\alpha=\psi\ast g_\alp, \qquad |\alpha|<k. \label{reproduce}
\end{equation}
\end{prop}

Note that for the expression $()^\gamma\ast'\psi=\sum_{j\in \Zd}
(\cdot-j)^\gamma\psi(j)$ to make sense $\psi$
needs to be continuous. The other two representations of $c(\gamma)$ are
valid for an arbitrary compactly supported distribution $\psi$.


\sl Proof. \rm  By Theorem \ref{repthm}, $\psi\ast'$ reproduces all
polynomials of degree $<k$, and hence (cf.\ e.g., \cite{BR})
$$\psi\ast' ()^\alp=\psi\ast ()^\alp=
()^{\alp}\ast'\psi,\quad \forall\alp\in\Zd_+,\ |\alp|<k.$$
Thus, $c(\alp)$ is well-defined.

Now, given $\alp$ as above, it is elementary that
(since $\psi\ast ()^0=\hatpsi(0)=1$)
$$\psi\ast ()^\alp=\sum_{|\bet|\le |\alp|}(\psi\ast ()^{\alp-\bet})(0)
\, ()^\bet=()^\alp+
\sum_{|\bet|< |\alp|}c(\alp-\bet) ()^\bet.$$
However, with $g_\alp$ as in (\ref{polrec}), we obtain (by convolving
$g_\alp$ with $\psi$, assuming by induction that $\psi\ast
g_\bet=()^\bet$ for $|\bet|<\alp$, and using the last identity) that
\begin{equation} \psi\ast g_\alpha\eqbd  \psi\ast ()^\alpha -
\sum_{|\beta|<|\alpha|} c(\alp-\bet)\,()^\bet=()^\alp.\end{equation}
\eop

Using this proposition with respect to the superfunction
$\psi\eqbd \sum_{\phi\in\Phi}\phi\ast'a_\phi$, we obtain the following:

\begin{thm} \label{reprFSI} Under the assumptions of Corollary \ref{fsipol},
the polynomials $(g_\alp)$ from Corollary \ref{fsipol1} satisfy the following
recurrence relation:
\begin{equation} g_\alpha\eqbd  ()^\alpha -
\sum_{\beta<\alpha} c(\alp-\bet)\,g_\beta,
\end{equation}
where
$$c(\gamma)\eqbd \sum_{\phi\in\Phi}c(\gamma,\phi),$$
while
$$c(\gamma,\phi)\eqbd  \sum_{j\in\Zd}(\phi\ast ()^\gamma)(j)a_\phi(-j).$$
Here, $\phi\ast ()^\gamma$ is {\sl continuous} convolution, while
$a_\phi$ is the finitely supported sequence that appears in Corollary
\ref{fsipol1}. Moreover, if each $\phi\in \Phi$ is continuous, we have
the alternative discrete convolution representation
$$c(\gamma,\phi)\eqbd  (\phi_{|_{\Zd}}^{\vphantom \gamma}
\ast a_\phi\ast ()^\gamma_{|_{\Zd}})(0)=
\sum_{j,k\in\Zd}\phi(j)\,a_\phi(k-j)\,(-k)^\gamma.$$
\end{thm}


\sl Remark. \rm Compare the last theorem with Theorem~1 of~\cite{CHM}.

\subsection{Good and bad superfunctions \label{goodbad1}}

Every FSI space contains a superfunction. This positive statement can be
turned negative: the existence of a superfunction in a given FSI space
tells us nothing about the structure of the space.  In contrast, 
Proposition~\ref{eq2propFSI}, Corollary~\ref{fsipol} and Theorem~\ref{reprFSI}
show that the existence of ``good'' superfunctions does lead us
to useful conclusions about the space and about the given generating set.
A particular useful condition is that the Fourier 
transform of the superfunction be bounded away from zero near the origin. 
This condition is important also from the numerical 
stability point of view. In view of the above, we say that a superfunction 
$\hatpsi$ is  {\bfi{nondegenerate}},\/ if  $\hatpsi$ near $0$ is bounded away 
from zero.

We note that the nondegeneracy itself falls short of classifying
``good'' superfunctions. For example, it can be checked that the 
$\l2$-projection of the sinc-function on an FSI space $S\subset \l2$ is 
always  nondegenerate as long as $S$ provides a positive approximation 
order~(by~\cite[Corollary 2.6]{BDR4}). However, the superfunctions obtained in 
this way may prove to be of little use due to their slow decay at $\infty$.
Thus, we require a complementary property of a superfunction: we say that
the superfunction has the {\bfi{finite span property}} if it is
in the {\it finite} span of the shifts of  the generating set
$\Phi$.  Such superfunctions are compactly  supported if $\Phi$ itself is.
It is proved in \cite{BDR2} that every local FSI space in $\l2$ contains a 
superfunction
that satisfies the finite span property (with $\Phi$ being any compactly
supported generating set for the space).

We call a superfunction {\bfi{good}},\/ if it is nondegenerate and finitely 
spanned by the shifts of $\Phi$. Such superfunctions are needed for 
constructing quasi-interpolation schemes. Indeed, the requirement appearing in 
Theorem~\ref{reprFSI} is exactly that the superfunction $\psi$ be good.

Corollary~\ref{sobcont} shows that any function in $S_\Phi(\ws)$ 
is of the form  $(v^*\widehat{\Phi})^\vee$, for  some $2\pi$-periodic 
vector-valued function
$v$.  Theorem~\ref{fsisob} adds that the vector  $v$ associated with 
a superfunction satisfies
$$ {v^* G_{\Phi,s}^0 v \over v^* G_{\Phi,s} v}=O(|\cdot|^{2k-2s}),$$
with $k$ the approximation order of the FSI space $S_\Phi(\ws)$.
The known $\l2$-theory of approximation orders of FSI spaces
offers then a recipe for constructing
good superfunctions: first solve the equation 
$G v=\widehat{\Phi}$  around the origin (cf. Theorem~\ref{conthm}) 
and then approximate $v$  by a trigonometric  polynomial vector $u$ such
that $v-u$ has a zero of order $k$ at the origin.
This is possible whenever 
the Gramian $G$ is $k$ times continuously differentiable around the 
origin and $G(0)$ is  invertible  (see~\cite[Theorem 4.2]{BDR4}).
Next, $v-u=O(|\cdot|^k)$ implies $v^*\widehat{\Phi} -u^*\widehat{\Phi}
=O(|\cdot|^k)$, hence $(u^* \widehat{\Phi})^\vee $ is a good superfunction.

Once $G(0)$ is not invertible, the notion of a ``good''
superfunction becomes more subtle. Are we only interested in the {\it
existence} of a superfunction
$\psi\in S_\Phi$ such that $\psi$ is ``reasonably local'' and
$\hatpsi(0)\not=0$, or do we also insist on simple ways to obtain that
function from the given generating set $\Phi$? Our discussion and development
focuses on the latter approach: after all, the SI space is {\it given} to us
in terms of the generating set $\Phi$, and we would like then the analysis
to stay as close as possible to this set.
Once we agree on that principle, it should be clear that ``very bad'' generating
sets $\Phi$ are not going to yield good superfunctions: for example
if $\Phi$ is compactly supported and $\hatPhi(0)=0$, there is no hope
to get from $\Phi$ in a simple way a superfunction $\psi$ with
$\hatpsi(0)\not=0$. The ultimate question is how to define ``good'' vectors
$\Phi$.  Our suggestion is simple: these are the vectors that yield
good superfunctions!

\bigskip
Our next results (in the next subsection) offer analysis
of vectors $\Phi$ whose Gramian is singular. We show the utility of this
analysis by providing a new proof to a famous example of de~Boor and H\"ollig
concerning the approximation order of $C^1$-cubics on three direction
mesh. We then provide an example of a ``seemingly good'' vector
$\Phi$ that cannot yield good superfunctions.

\subsection{\label{estim}Estimating approximation orders when  Gramians are 
singular}

Theorem~\ref{fsisob} enables us, at least in principle, to determine the
order of approximation provided by a given (stationary or nonstationary)
FSI ladder in any Sobolev space. But, as we already saw in the $\l2$-case, 
such analysis is hard to carry out if the Gramian $|\cdot|^{-2s}G_{\Phi,s}$ 
is singular at the origin. The problem is exacerbated by the fact that the 
entries of $G_{\Phi,s}$ may be hard to compute.

Let us examine closely the source of the difficulty. For $s=0$, the
computation of 
approximation orders depends on estimating ratios of the form
$${v^\ast G^0 v\over 
v^\ast Gv}$$
around the origin. Without loss, one can assume that the vector
$v$ is normalized pointwise. If $G$ is continuous at $0$ and invertible there,
we can then dismiss the denominator, since it does not affect the asymptotic
behaviour of the above expression. In contrast, if $G$ is singular at the
origin, the denominator might affect the
approximation order. The use of the verb ``might'' is justified: roughly
speaking, there is hope that the specific vectors $v$ that minimize the 
numerator are far enough from the kernel of $G(0)$. Whenever this is the
case, the problem is reduced to examining the behaviour of the numerator only.
The current subsection translates the above discussion into  rigorous
analysis.

We first provide below a theorem that establishes an upper bound on 
the approximation order of an FSI space. The upper bound does not require the
invertibility of the associated Gramian.
To this end, we denote by $\sigm_{\mi{}} (A)$ the smallest eigenvalue
of a positive-definite Hermitian matrix $A$.\index{$\sigm_{\mi{}} (A)$ --
the smallest eigenvalue of a Hermitian matrix (function) $A$}

\begin{thm} \label{estthm} Suppose $\widehat{\Phi}\subset 
L_\infty(\Omega)$  for some neighborhood $\Omega$ of the origin.
Given any set ${\cal I}\subseteq \nzero$, denote by $k(\Phi,{\cal I},s)$ 
the order of the zero that the scalar function
$$\omega\mapsto
\sigm_{\mi{}}(A(\ome)),\quad A(\ome):=
\sum_{\alpha \in {\cal I}} \widehat{\Phi}(\omega +\alpha) 
\widehat{\Phi}(\omega +\alpha)^* |\omega+\alpha|^{2s}$$ 
has at the origin.
Then the approximation order provided by
the FSI space $S_\Phi(\ws)$ is no larger than $k(\Phi,{\cal I},s)/2$.
\end{thm}

\sl Proof. \rm Suppose that $S_\Phi$ provides approximation order
$k$ in $\ws$. Then the characterization from Theorem~\ref{fsisob} 
implies that, for 
$v_0:\omega \mapsto v_0(\omega)$ that 
minimizes~(\ref{lstat}) pointwise, the expression
$$ {v_0^* G^0_{\Phi,s} v_0^{\vphantom *} 
\over v_0^* G_{\Phi,s} v_0^{\vphantom *} } {1\over |\cdot|^{2k-2s}}=:\Mu_{\Phi,s}$$
is bounded in  a neighborhood $\Omega$ of the origin. Using the identity
$$ (1-\Mu_{\Phi,s} |\cdot|^{2k-2s} ) v_0^* G^0_{\Phi,s} v_0^{\vphantom *} 
=\Mu_{\Phi,s}  |v_0^*\widehat{\Phi}|^2 \, |\cdot|^{2k},$$ 
and the fact that $k>s$, we conclude that  
$v_0^* G_{\Phi,s} v_0^{\vphantom *}$
is bounded above by a constant multiple of $|v_0^*\widehat{\Phi}|^2 \,
 |\cdot|^{2s}$ (compare with the argument used in the proof of
 Proposition~\ref{eqprop}). Thus, 
$$  \inf_v {v^* G^0_{\Phi,s} v 
\over v^* G_{\Phi,s} v } =  {v_0^* G^0_{\Phi,s} v_0^{\vphantom *} 
\over v_0^* G_{\Phi,s} v_0^{\vphantom *} } \geq \const {v_0^* G^0_{\Phi,s}
v_0^{\vphantom *} 
\over |v_0^* \widehat{\Phi}|^2 |\cdot|^{2s}}\geq
\const  \inf_{|v|=1} {v^* G^0_{\Phi,s} v 
\over |v^* \widehat{\Phi}|^2 |\cdot|^{2s}} \geq \const |\cdot|^{-2s} 
\inf_{|v|=1} v^* G_{\Phi,s}^0 v.$$ The last inequality uses the 
assumption  that $\widehat{\Phi}$ is 
bounded around the origin. We conclude then that, for some $C>0$ and
a.e. around the origin,
$$  \inf_v {v^* G^0_{\Phi,s} v \over v^* G_{\Phi,s} v } 
\geq C |\cdot|^{-2s} \sigm_{\mi{}} (G^0_{\Phi,s}).$$
But 
$G_{\Phi,s}^0 -  A$
is  (pointwise) a nonnegative definite Hermitian matrix, hence, pointwise,
$\sigm_{\mi{}}(G^0_{\Phi,s}) \geq \sigm_{\mi{}}(A).$
The desired result then follows from  Theorem~\ref{fsisob}.
\eop

As alluded to before, we know quite precisely when the above upper bound
matches the associated approximation order.

\begin{thm}  \label{estthm2} Suppose $\widehat{\Phi}\subset 
L_\infty(\Omega)$  for some neighborhood $\Omega$ of the origin. 
Let $v_0$ be  a normalized eigenvector of $G_{\Phi,s}^0$
associated with the minimal eigenvalue (i.e., for a.e. $\ome\in\Ome$,
the pair $(\sigm_{\mi{}}(G_{\Phi,s}^0(\ome)),v_0(\ome))$ is an eigenpair
of $G_{\Phi,s}^0(\ome)${\rm )}. If
$|v_0^*\widehat{\Phi}|$ 
is bounded away from zero almost everywhere in $\Omega$, then the
approximation order of $S_\Phi(\ws)$ is exactly $k(\Phi,\nzero,s)/2$.
\end{thm}

\sl Proof. \rm One only needs to show that the approximation order
of $S_\Phi$ is bounded below by $k(\Phi,\nzero,s)/2$.
But $$\inf_v {v^* G^0_{\Phi,s} v 
\over v^* G_{\Phi,s} v } \leq  {v_0^* G^0_{\Phi,s} v_0 
\over v_0^* G_{\Phi,s} v_0 } =  {\sigm_{\mi{}} (G^0_{\Phi,s}) 
\over |v_0^* \widehat{\Phi}|^2 |\cdot|^{2s} + \sigm_{\mi{}} (G^0_{\Phi,s})} 
\leq {\sigm_{\mi{}} (G^0_{\Phi,s}) 
\over |v_0^* \widehat{\Phi}|^2 |\cdot|^{2s}} \leq  
\const |\cdot|^{-2s}\sigm_{\mi{}} (G^0_{\Phi,s}) .$$
Theorem~\ref{fsisob} then yields the requisite lower bound.
\eop

\subsection{\label{c1cubic}Example: bivariate $C^1$-cubics}

The results of the last section raise two questions. The first is whether
the upper bounds provided in Theorem~\ref{estthm} are useful, i.e.,
whether they can be applied to solve a nontrivial problem. We provide
in the current subsection an affirmative answer to this question.

The other, more fundamental question, is whether the setup of
Theorem~\ref{estthm2} is universal, i.e., whether we can {\it always} 
dispense with the denominator in the characterization provided in
Theorem~\ref{fsisob}. This question is intimately related to the existence
of good superfunctions. In the next subsection we provide a (-n unfortunate)
negative answer to that second question.

As said, we describe now
an example where Theorem~\ref{estthm} applies in a 
situation when direct evaluation of the approximation order is quite 
complicated.  We choose the notorious example of an FSI space that reproduces 
all polynomials of order $\le 3$, but provides only approximation order $3$. The
example first appears in~\cite{BH}. A second, completely different, proof
of this result appears in \cite{BDR4}. Our proof is thus the third one for
this result.

Consider  the following two bivariate
compactly supported piecewise polynomial functions whose
Fourier transforms are given by
\begin{eqnarray*}
&\widehat{\phi}_1(u,v)= \ii{(v(1-e^{-\ii w})-w(1-e^{-\ii v}))(1-e^{-\ii u})
(1-e^{-\ii v})(1-e^{-\ii w}) \over  (uvw)^2}, &\\
&\widehat{\phi}_2(u,v)=  \widehat{\phi}_1(v,u),&
\end{eqnarray*} 
where $w=u+v$.  These functions are known as the Fredrickson 
elements. With $\Phi\subset \l2(\R^2)$ the 2-vector consisting of the 
above functions, the Gramian $G_\Phi$ is invertible in a punctured
neighborhood of the origin. Hence it is still possible to compute enough 
coefficients in the Taylor expansion of $1-\widehat{\Phi}^* G^{-1}_\Phi 
\widehat{\Phi}$
to find the first nonvanishing nonconstant term. This complicated
analysis was carried out in \cite{BDR4}. It  shows
that the first nonzero term in that expansion if
of order $6$, so $S_\Phi$ provides approximation order $3$.

We use here, instead, Theorem~\ref{estthm} to arrive at the
same conclusion more easily. Let ${\cal I}=\{(0,2\pi),(2\pi,0) \}$. 
Then 
\begin{equation}
\sum_{\alpha \in {\cal I}} \widehat{\Phi}((u,v) +\alpha)
\widehat{\Phi}((u,v)
 +\alpha)^*= \Psi_1(u,v) \Psi_1^*(u,v)+\Psi_2(u,v) \Psi_2^*(u,v)+
 o((|u|^2+|v|^2)^3)\label{here}
 \end{equation}
 where
 \begin{eqnarray*}
 \Psi_1(u,v)={(1-e^{-\ii u})(1-e^{-\ii v})(1-e^{-\ii(u+v)})\over
 (u+v+2\pi)^2} \left[ \begin{array}{c} {-2\pi+\pi \ii u +2 \pi \ii v
 +\ii uv/2+\ii v^2/2 \over  u (v+2\pi)^2 } \\
 { 2\pi-\pi \ii u +\ii u v/2+\ii v^2/2 \over u(v+2\pi)^2}
 \end{array} \right] , \\
 \Psi_2(u,v)={(1-e^{-\ii u})(1-e^{-\ii v})(1-e^{-\ii(u+v)})\over
 (u+v+2\pi)^2} \left[ \begin{array}{c}
 { 2\pi-\pi \ii v +\ii uv/2+\ii u^2/2 \over v(u+2\pi)^2} \\
 {-2\pi+\pi \ii v +2 \pi \ii u +\ii uv/2 +\ii u^2/2 \over v (u+2\pi)^2 }
 \end{array} \right].
 \end{eqnarray*}
The trace of the matrix $\Psi_1^{\vphantom *} \Psi_1^*+\Psi_2^{\vphantom *} \Psi_2^*$ is of order $4$, 
whereas its determinant is of order $10$, so its minimal eigenvalue vanishes
to order $6$ and its maximal eigenvalue to order $4$.  Since the (matrix) terms
that were left out of the computation are all of order $o(|\cdot|^6)$, the 
eigenvalues of the left hand side of~(\ref{here}) are also of order
$6$ and $4$, respectively. Now Theorem~\ref{estthm} implies that the 
approximation order of $S_\Phi$ is at most $3$.  

The fact that the approximation order is at least $3$ is trivial:
the sum $\phi_1+\phi_2$ yields a superfunction which is nothing
but the box spline $M_{2,2,1}$ (whose approximation
order is indeed $3$). The vector $\Phi$ is thus an example where the 
singularity of the Gramian does not preclude the existence of a good 
superfunction.

From our standpoint, the $C^1$-cubic vector $\Phi$ is
``good'', since the finite span of its shifts contains  a good superfunction.
The notoriety of this case is due to the difficulty in asserting
that this $\Phi$ provides approximation order no higher than $3$.
The fact that the space reproduces all cubic polynomials is a sad,
misleading, accident. The reader may claim that we ignore the fact
that $\Phi$ here provides an approximation order which is ``disappointing''.
While that might be the case, it goes beyond the realm of this article:
we are only interested in ways to capture the approximation order of the
given space, and not in the construction of SI spaces that
provide ``satisfactory'' approximation order.

\subsection{Good and bad superfunctions, continued \label{goodbad2}}

We will now show an example of a vector $\Phi$ whose entries
seem to be ``reasonable'' but nonetheless it does not admit a good
superfunction. This example of a bad $\Phi$,
together with the example from the last
section of a good $\Phi$, illustrates the depth of the difficulty
in pinning down the notion of a good generating set for an FSI space.

Let $g$ be a compactly supported bivariate function whose Fourier
transform $\hatg$ has a zero of order $k>2$ at each of the $\nzero$
points. Moreover, we assume that $1-\hatg$ has a zero of order $k$
at the origin. There are many ways to construct such a function.
For example, one can take the univariate B-spline of order $k$,
apply a suitable differential operator $p(D)$ or order $k-1$ to it, and then
use its tensor product in 2 dimensions. In this case,
$g$ is piecewise polynomial of local degree $k-1$ in each of its
variables, and with support $[0,k]^2$. It provides approximation order
$k$ (in $L_2$, for example).

Now, let $e$ be the bivariate exponential with
frequency $(2\pi,0)$, i.e., $e:x\mapsto e^{2\pi \ii x(1)}$.
We define a vector $\Phi$ with two components
$$\phi_1\eqbd g+e D^{(0,2)}g,\quad \phi_2\eqbd g-e D^{(2,0)}g.$$
Here, to recall, $D^\alp$ is the normalized monomial differentiation,
viz., $2D^{(2,0)}g$ is the second derivative of $g$ in the first argument.
Despite the fact that each $\phi_i$ provides only approximation
order $2$, we contend that the FSI space $S_\Phi(\l2)$ provides approximation 
order $k$.  We construct, to this end, a compactly supported superfunction as
follows.

We choose a vector-valued function $v$ with two components that are 
trigonometric polynomials such that 
$$v-\pmatrix{()^{2,0}\cr ()^{0,2}\cr}=O(|\cdot|^{k+2})$$
around the origin. We then note that the Taylor expansion of order $k$
of $\hatPhi$ around the point $(2\pi,0)$ is
$${1\over 2}\pmatrix{-()^{0,2}\cr {\hphantom {-}}()^{2,0}\cr}.$$
At the same time, the Taylor expansion of order $k$ of $\hatPhi$ around
any point of $\nzero$ other than $(2\pi,0)$ is zero. From this, we conclude
that the compactly supported function $\psi$ that is defined by
$$\hatpsi\eqbd v^\ast\hatPhi$$
has a zero of order $k+2$ at {\it every} point of $\nzero$. Finally,
at the origin, $\hatpsi-|\cdot|^2\hatg$ has a zero of order $k+2$.

In order to determine the approximation order provided by $\psi$, we consider
the expression
$${[\hatpsi,\hatpsi]-|\hatpsi|^2\over |\hatpsi|^2}=
{[\hatpsi,\hatpsi]-|\hatpsi|^2\over |\cdot|^4|\hatg|^2}
{|\cdot|^4|\hatg|^2\over |\hatpsi|^2}.$$
The term
$${|\cdot|^4|\hatg|^2\over |\hatpsi|^2}$$
is bounded around the origin. The other term,
$${[\hatpsi,\hatpsi]-|\hatpsi|^2\over |\cdot|^4|\hatg|^2}$$
remains bounded (around the origin) even when multiplied by $|\cdot|^{-2k}$.
Thus, $\psi$ provides approximation order $k$ (in $L_2$), {\it a fortiori}
$\Phi$ provides that approximation order.

Note that the superfunction $\psi$ does not satisfy the desired
condition $\hatpsi(0)\not=0$. In fact, this is necessary in a certain
sense. Indeed, let $\tau$ be  a $2\pi$-periodic vector-valued function 
that is continuous at the origin and does not vanish there.
Let us further define a  function $f$ by
$$\hatf=\tau^\ast\hatPhi.$$
Then, up to a non zero multiplicative constant,
the low order derivatives of $\hatf$ at $(2\pi,0)$ coincide with
the derivatives at the origin of the function
$$(\tau_1 ()^{(0,2)}-\tau_2 ()^{(2,0)})\hatg.$$
Since we assume $\tau$ not to vanish at the origin, it is clear that some
second order derivative of the above expression does not vanish at the origin.
As such, $f$ cannot provide approximation order larger than $2$.

While the vector $\Phi$ in this example does not yield a good superfunction,
it satisfies the following positive property: we could use the truncated
Gramian $G_{\Phi,0}^0$ is order to determine the approximation order
of $S_\Phi$. Indeed, if we normalize the given vector $v$
(i.e., redefine it pointwise as $v/|v|$), we obtain a vector for which
$v^\ast G_{\Phi,0}^0v$ yields the correct decay rate ($k$) at the origin.
This means, in turn, that the smallest eigenvalue of $G^0_{\Phi,0}$
still determines the approximation order of the space $S_\Phi$.
The superfunction that we obtain in this way (i.e., by using the normalized
$v$) is still not good: it decays painfully slowly at $\infty$.

\medskip
We close this section with two comments:
\begin{itemize}
\item We do not know at present of an example where the smallest eigenvalue
of the truncated Gramian $G_{\Phi,s}^0$ does not determine the approximation
order of the space provided, of course, that $\hatPhi(0)\neq 0$.
\item The above example (i.e., of a case when the smallest eigenvalue
of $G_{\Phi,s}^0$ determines the approximation order while there exists no
good superfunction) is very much a multivariate phenomenon. It is not hard
to prove that such a case is impossible in one variable, and we leave it as an
exercise to the interested reader.
\end{itemize}

\subsection{An application: approximation orders of smooth refinable 
functions \label{smoothsec1}}

We provide in this section one of the most interesting applications 
of superfunction theory: lower bounds on approximation orders of
smooth refinable vectors. We note that approximation orders of
refinable vectors are treated in more detail in the next section. 
However, the current topic fits better into the realm of this section.

At base, our result will show that once $\Phi$ is refinable, and once $S_\Phi$
contains a single nonzero function $\psi$ from a certain class,
the stationary ladder generated by $\Phi$ must provide an approximation 
order that corresponds to the class of $\psi$. Our definition of the
``class'' in question requires the Fourier transform of $\psi$ to decay
(in a weak sense) at a certain rate.

This problem has rich history in the context of PSI ladders
(see the introduction to~\cite{R}).  A substantial treatment of the FSI 
case is given in~\cite{R}. However, that treatment is carried out under
the assumption that the Gramian $G_\Phi$ is invertible at the origin.
In contrast, we focus in this paper on the situation where there are
multiple solutions to a single refinement equation, and in such a case the
Gramian of any particular solution is {\it not} invertible at the
origin. This understanding was our motivation to look for an alternative
approach to that of \cite{R}. It is useful to stress that, in general, 
refinable vectors that contain a smooth (even analytic!) function need 
not provide any positive approximation order at all.  (An example of this 
type can be found in \cite{R}.) Thus, one must impose certain side conditions
either on the vector $\Phi$ or on the function $\psi$.

\medskip
Let $P$ be an $r\times r$ matrix whose entries are $2\pi$-periodic and 
measurable. Let $\Phi$ be a vector-valued function with $r$ components
whose entries are in $\ws$ for some $s\in\R$. We say that $\Phi$ is 
{\bf refinable} if the functional equation
\begin{equation}
\hatPhi(2\cdot)=P\hatPhi\label{balagan}
\end{equation}
is satisfied.

Our goal is to prove the following result.
In the result, as well as elsewhere in this subsection, we use the
following notation:
\begin{equation}A\eqbd \{\ome\in\Rd:\ {1\over 2}\rho< |\ome|\le \rho\}.
\label{defA}\end{equation}
Here, $\rho\in (0,\pi)$ is arbitrary, but fixed. 

\begin{thm} \label{thmref1} Let $s\le 0$, and let
$\Phi\subset \ws$ be a solution to~(\ref{balagan}). With $A$ as in
(\ref{defA}), assume that there exists $f\in S_\Phi(\ws)$
with the following properties:
\begin{enumerate}
\item $|\hatf|$ is bounded above as well as away from zero on $A$.
\item The numbers
$$\lam_m\eqbd \norm{\sum_{\alp\in 2^m( 2\pi\sZZ^d\bks0)}
|\hatf(\cdot+\alp)|^2}_{L_\infty(A)}, \quad m\in\Z_+$$
satisfy  $\lam_m=O(2^{-2mk})$, for some positive $k$.
\end{enumerate}
Then $S_\Phi(\ws)$ provides approximation order $k$.
\end{thm}

We approach this result via the notion of
the {\bfi{dual equation\/}}  \begin{equation}
v^*(2\cdot)P= v^*.\label{dual}
\end{equation}
Here, $v$ is a vector-valued function with $r$ components.
We require the dual equation to be valid in some (small) ball  $U$ 
centered at the origin, and define the entries of the {\bf dual vector} $v$ 
to be equal to $0$ outside $[-\pi,\pi]^d\bks U$. We then extend $v$ to a 
$2\pi$-periodic vector. Thus, $v$ is supported on $U+2\pi\Zd$, and 
satisfies~(\ref{dual}) there.

We collect in the next lemma a few simple facts about dual vectors.

\begin{lem} \label{observe} Let $A$ be as in (\ref{defA}).
Given any $v_0$ defined on $A$,
equation~(\ref{dual}) can be solved on the punctured
disk 
\begin{equation}U\eqbd  \{ \omega : 0< |\omega|\leq  \rho \}
\label{defD}
\end{equation} so that the solution $v$ 
satisfies $v|_A=v_0$.  Moreover, we have then,
a.e.\ on $U$,
$$ v^*(\ome/2^m) \widehat{\Phi}(\ome/2^m+\alpha) =v^*(\ome) 
\widehat{\Phi}(\ome + 2^m\alpha), \quad {\rm all} \; m\in \Z_+, \;\; \alpha 
\in 2\pi\Z^d.$$
\end{lem}

\sl Proof. \rm  We define $v$ by
 $v^*(\omega)\eqbd 
v^*(2\omega) P(\omega)$, for all $\omega \in 2^l A$,
$l=-1, -2, \ldots$. Then $v$ clearly satisfies (\ref{dual}) (on $U$, and hence
on $U+2\pi\Zd$).
 
The second part of the lemma is obtained by iterating $m$ times
with
\begin{eqnarray*} && v^* \widehat{\Phi} (\cdot+ \alpha)= 
v^*(2\cdot)P \widehat{\Phi} (\cdot+\alpha)= 
v^*(2\cdot)P(\cdot +\alpha) \widehat{\Phi} (\cdot+\alpha)= \\
&&  v^* (2\cdot) \widehat{\Phi}(2
 (\cdot+\alpha)). 
\end{eqnarray*}  \eop

\goodbreak

\medskip \sl Proof of Theorem \ref{thmref1}. \rm 
By Corollary~\ref{sobcont},
$\hatf=\tau^\ast\hatPhi$ for some $2\pi$-periodic $\tau$.
Denoting by $v_0$ the restriction of $\tau$ to $A$, we extend $v_0$
to a dual vector $v$ by Lemma~\ref{observe}.   Defining $\psi$
by
$$\hatpsi\eqbd v^\ast\hatPhi,$$
we have, by the same lemma, that,
for a.e. on $A$, and for every nonnegative integer $m$,
$$|\hatpsi(\ome/2^m)|=|\hatf(\ome)|.$$ 
Thus, in view of our assumptions on $f$, we conclude that
$|\hatpsi|$ is bounded between two positive constants around the origin.

Next, we prove that
$S_\psi(L_2)$ provides approximation order $k$.  The argument will show, as a by-product,
that $\psi\in \ltwo$. 

Since $\hatpsi$ is bounded away from $0$ around the origin, it remains to 
prove, in view of Result \ref{psil2}, that 
$[\hatpsi,\hatpsi]^0|\cdot|^{-2k}$ is bounded around the origin
(with $[\hatpsi,\hatpsi]^0\eqbd [\hatpsi,\hatpsi]-|\hatpsi|^2$).
Let $\ome\in A$, and $m$ a positive integer. Then, by the definition
of $\psi$ and Lemma \ref{observe}
$$[\hatpsi,\hatpsi]^0(\ome/2^m)=\sum_{\alp\in 2\pi\Zd\bks0}
|\hatpsi|^2(\ome/2^m+\alp) =\sum_{\alp\in 2\pi\Zd\bks0}
|\hatf|^2(\ome+2^m\alp)\le \lam_m\le C 2^{-2mk}\le C |\ome/2^m|^{2k}.$$
Thus, $[\hatpsi,\hatpsi]^0=O(|\cdot|^{2k})$, on the punctured
ball $U$ of radius
$\rho$ centered at the origin. Since $\hatpsi$ is supported on
$U+2\pi\Zd$, it follows that $\hatpsi$, hence $\psi$, lies in $\ltwo$.
Result \ref{psil2}
then applies to show that $\psi$ provides approximation order $k$ in
$\ltwo$.

On the other hand, $\hatpsi=v^\ast\hatPhi$, with $v$
measurable and $2\pi$-periodic. Since $\psi\in\ltwo$, and
$s\le 0$,  we have that $\psi\in \ws$. 
Corollary~\ref{sobcont} shows then  that $\psi\in S_\Phi(\ws)$.
Now, Proposition~\ref{eqprop} implies that $\psi$ provides approximation
order $k$ in $\ws$, hence $\Phi$ provides also approximation order $k$ in
that space.
\eop

\begin{dis} \label{setA} 
As the proof of the theorem shows, there is in fact more
freedom in the choice of $A$. It suffices to assume that $A$ is compact,
that the intersection $A\cap 2A$ has measure zero, and that the union
$\cup_{m=-\infty}^0 A/2^m$ contains a neighborhood of the origin. The proof 
remains essentially the same.  
\end{dis}

\section{Vector refinement equations \label{refinesec}}

In our studies so far, we considered SI spaces one at a time.
There are situations, however, where several different SI spaces may stem 
from one common source. In cases of this type, it is important to study
the resulting SI spaces in a cohesive, combined, way.

The best examples of this type are the multiple vector-valued solutions 
to refinement equations, and this  is, indeed, the topic of the current 
section.  Let us start with the requisite definitions.

Let $P$ be an $r\times r$ square matrix whose entries are,
$2\pi$-periodic (measurable) functions (defined on $\Rd$). The 
functional equation 
\begin{equation}
\hatPhi(2\cdot)=P \hatPhi, 
\label{*}
\end{equation}
is a {\bfi{vector  refinement equation}}, $P$\index{$P$ -- refinement mask} 
is a {\bfi{refinement (matrix) 
mask}}, and a solution $\Phi$ is a {\bfi{refinable vector}}.
Here, the entries  of the vector $\Phi$ are (measurable) functions, 
or, more generally, tempered distributions, defined on $\Rd$. The rows 
and columns of the matrix $P$ are, thus, indexed by either
the integers $1,\ldots,r$, or, more directly, by the entries of $\Phi$.
In this generality, the equation~(\ref{*}) has, as a rule, infinitely
many linearly independent solutions. Indeed, if $P$ is nonsingular around 
zero, then a solution $\hatPhi$ can be chosen arbitrarily on a set $A$
of the `dyadic annulus' type introduced in Discussion~\ref{setA} and then 
continued to the rest of the Fourier domain using the recipe 
\begin{eqnarray*} \hatPhi(2\ome)  \eqbd  P(\ome)
\hatPhi(\ome), &  \ome\in 2^jA,  &  j=0,1,2,\ldots \\
  \hatPhi(\ome)  \eqbd   P^{-1}(\ome)\hatPhi(2 \ome), & \ome\in 2^j A, 
& j=-1,-2,\ldots. 
\end{eqnarray*}
Most of the solutions of the above type will decay very slowly (will not be
even in $L_1(\Rd)$). In contrast, if we assume the entries of
$P$ to consist of trigonometric polynomials, and if we correspondingly
insist on compactly supported  solutions, then the solution
space is necessarily finite-dimensional (as explained in detail in the next 
section). The special instance when the compactly supported solution
space is one-dimensional is quite well understood (see, 
e.g.,~\cite{BDR4},~\cite{CHM1},~\cite{J1}, and~\cite{JS}).
We are therefore primarily interested in the case when there are multiple 
(in a nontrivial sense) compactly supported solutions  to the 
equation~(\ref{*}). We denote by $$R(P)$$\index{$R(P)$ -- 
compactly supported solutions to the refinement equation with mask $P$} 
the linear space of all the 
solutions of (\ref{*}) whose entries are compactly supported distributions. 

The core of study here is the connection between properties
of the refinement mask $P$ and its corresponding solution(s) $\Phi$.
Shift-invariant spaces enter the discussion in a very natural way. For
example, if the solution vector $\Phi$ lies in $L_2$, then
one has the inclusions
$$S_\Phi(L_2)\subset \calD(S_\Phi(L_2)),$$
with $\calD$ the dilation operator $f\mapsto f(2\cdot)$.
Due to the above inclusion, we refer to the SI spaces $S_\Phi$ generated
by a refinable $\Phi$ as a {\it refinable SI space}.

 
We start our study in this section with the problem of
{\it{existence}} of compactly supported solutions to~(\ref{*}). 
Our second, and main, topic is the characterization of the
approximation orders of the FSI space generated by the solutions $\Phi$
to the refinement equation. This study is based on the premise that,
in the case where multiple solutions to the same equation exist,
the objective should be the interplay among those solutions, and not only
the individual properties of each one of them. In this course of study,
we introduce the notions of the {\it{combined Gramian}} and
the {\it{coherent approximation order}} and connect them
with the (i) the
approximation orders of the SI spaces generated by the solutions to the
equations, (ii) the polynomial reproduction property of the mask $P$, and
(iii) the sum rules satisfied by $P$. Finally, we already
provided (in \S\ref{smoothsec1}) lower bounds on the approximation order of 
a refinable SI space in terms of the smoothness of the smoothest function
in that space. 

\def\zn{{\cal Z}_N}
\def\zm{{\cal Z}_M}
\def\alp{\alpha}
\def\bet{\beta}

\subsection{\label{rfsec}Compactly supported solutions to the refinement
equation}

The structure of the compactly supported solutions of~(\ref{*}) was first 
completely described in~\cite{JJS}. We now restate the main result of that 
paper and provide a different proof for it. We use the partial
order $\le $ on $\Z^d$, defined by\index{$\le$ -- partial order on $\Z^d$} 
$$ a\le b  \Longleftrightarrow a-b\in \Z_+^d.$$
Also, given a nonnegative integer $N$, we set\index{$\zn$ -- multi-indices of length smaller than $N$} 
$$\zn\eqbd  \{\alp\in \Z^d_+:\ |\alp|\eqbd \alp_1+\ldots+\alp_d\le N\}.$$
Finally, we recall that the definition of the monomial differential operator 
$D^\gamma$ includes the normalization factor $1/\gamma!$.

\begin{thm} \label{charthm}
Given  an $r\times r$-matrix $P$ whose entries are 
trigonometric polynomials, set 
$$N\eqbd  \max \{ n : 2^n \in \spec(P(0))\}.$$
Then the map  
$$ \Phi \mapsto ((D^\alpha\widehat{\Phi})(0))_{\alpha\in\zn},$$
is a bijection between the collection $R(P)$ of all compactly supported 
solutions of~(\ref{*}) and the kernel $\ker L$ of the map
$$L : \C^r\times \zn \to  \C^r\times \zn
: (w_\alpha)_{\alp\in\zn} \mapsto (2^{|\alpha|}w_\alpha-
\sum_{0\leq \beta \leq \alpha} 
(D^{\alpha - \beta}P)(0)\,w_\beta:\  \alpha\in \zn).$$
\end{thm}

\sl Proof. \rm Let $\Phi$ be a compactly supported distributional 
solution to~(\ref{*}), and denote $w_\alp\eqbd (D^\alp\hatPhi)(0)$,
$\alp\in\Z_+^d$. Since the vector-valued function $\widehat{\Phi}$ 
is entire, the vectors $w_\alp$ are all well-defined. Moreover, one 
easily concludes 
from the relation~(\ref{*}) (by applying $D^\alp$ to both sides
of that identity, expanding the right hand side with the aid
of Leibniz' formula, and evaluating the result at
$0$) that the sequence  $(w_\alp)_{\alp\in\Z_+^d}$ solves
the infinite triangular system 
\begin{equation}
2^{|\alpha|}w_\alpha=\sum_{0\leq \beta \leq \alpha}
 (D^{\alpha -\beta}P)(0) w_\beta, \qquad \alpha \in \Z_+^d.   \label{**}
\end{equation}
In particular, $(w_\alp)_{\alpha\in\zn} \in \ker L$.

Conversely, let $w\eqbd  (w_\alpha)_{\alpha\in\zn}\in \ker L$. Then $w$ 
extends uniquely to a solution to~(\ref{**}) (in order to solve uniquely
for $w_\alp$ in~(\ref{**}) one needs the matrix $2^{|\alp|}I-P(0)$ to be
invertible, which is indeed the case for every $|\alp|>N$, by our assumption
of $\spec P(0)$).

Let $\| \cdot \|$ be any vector norm on $\C^r$. The operator norm on
$\C^{r\times r}$ subordinate to $\| \cdot\|$ will be denoted in the same way.
We claim that for some constant $A>0$, 
\begin{equation} \|w_\alpha\|\leq A^{|\alpha|}/\alp!, \qquad {\rm all} \quad
\alpha\in \Z_+^d. \label{***} \end{equation}  
Let us see first that (\ref{***}) yields the existence of a suitable solution
to (\ref{*}).

With (\ref{***}) in hand, we define (with $()^\alp$ the normalized monomial)
$g\eqbd  \sum_{\alpha \in \Z_+^d} \alpha! ()^\alpha\,w_\alpha$, and observe
that (each of the entries of)
$g$ is entire of exponential type, i.e., it satisfies
\begin{equation}
\|g(\omega)\| \leq \ee^{\widetilde{A}|\omega|}, \qquad {\rm all} \quad \omega \in \C^d,
\label{****}
\end{equation}
where $|\cdot |$ denotes an arbitrary norm on $\C^d$. We need further to show
that each of the entries of $g$ is the Fourier transform of a compactly
supported distribution, which, by the Paley-Wiener-Schwartz Theorem
(Theorem~19.3 from~\cite[p.~375]{Ru} and Exercise~7.4
from~\cite[p.~216]{Ra})
amounts to showing that (in addition to (\ref{***})) the restriction of
$g$ to $\Rd$ has slow growth at $\infty$.  
In order to prove the requisite slow growth,
we follow an argument from~\cite{JS}:  Denoting
$$C_1\eqbd  \sup_{\xi \in \R^d} \|P(\xi)\| \quad \hbox{and}\quad C_2\eqbd  
\sup_{1\leq |\xi| \leq 2} \|g(\xi)\|,$$
we pick $\omega \in \R^d$,  such that $1\le | \omega |< 2$. 
By the construction of $g$, $g(2\cdot)=Pg$, and hence, for every positive
$n$, $g(2^n\ome)=P({2^{n-1}\omega}) \cdots 
P({\omega }) g({\omega })$. Consequently,
$$ \| g(2^n\omega)\| \leq C_2 C_1^n \le  C_2 (2^n|\omega|)^{\log_2(C_1)},$$ 
a bound that evidently establishes the sought-for slow growth.

It remains to prove (\ref{***}). To that end,
we pick $N_0\in \N$ so that $$2^{-N_0} \|P(0)\|<1, \qquad \left({3\over 4} 
\right)^{N_0} {1 \over 1-2^{-N_0} \|P(0)\|}\leq 1.$$ 
Since $P$ is a matrix of trigonometric polynomials, there exists
$A>0$ such that $\| (D^\alpha P)(0)\|\leq (A/2)^{|\alpha|}/\alp!$, all 
$\alpha$. Moreover, by modifying $A$ if need be,
we may assume that $A$ satisfies the estimate~(\ref{***})
for every $\alp\in {\cal Z}_{N_0}$. 

In order to prove~(\ref{***}) for $|\alp|>N_0$, we may assume,
by induction, that~(\ref{***}) holds for all $\beta$ 
with $|\beta|<|\alp|$. Then, by~(\ref{**}),
$$ \|(2^{|\alpha|}I-P(0))w_\alpha \|\leq \sum_{0\leq \beta < \alpha} 
 {A^{|\alpha -\beta|}\over 2^{|\alp-\bet|} (\alp-\bet)!} 
 {A^{|\beta|}\over \bet!} = A^{|\alp|} \sum_{|\bet|<|\alp|}
 {1\over  2^{|\alp-\bet|} (\alp-\bet)!  \bet!} 
 \leq {A^{|\alpha|}\over \alp!} \left( {3 \over 2}\right)^{|\alpha|}, $$
hence
$$ \|w_\alpha\|\leq \left( {3 \over 4}\right)^{|\alpha|} 
{A^{|\alpha|} \over \alp! (1-2^{-|\alpha|} \|P(0)\|)} \leq
{A^{|\alpha|} \over \alp!}.$$
This proves~(\ref{***}), and the proof is thus complete.   \eop

The theorem can be extended to 
refinement equations more general  than (\ref{*}). For example, we can replace the 
dilation by $2$ by a dilation by any matrix $M$ which is {\it{expansive}}, i.e., 
its spectrum lies outside the closed unit disc.

\begin{thm} \label{char2thm} 
Given  an $r\times r$-matrix $P$ whose entries are 
trigonometric polynomials and an expansive $d\times d$ matrix $M$, set 
$$N\eqbd  \max \{ n : 0\in \spec\{M^n-P(0)\}\}.$$
Then, the map 
$$ \Phi \mapsto ((D^\alpha\widehat{\Phi})(0))_{\alpha\in\zn},$$
is a bijection between the collection of all (compactly supported) solutions
of the refinement equation $\hatPhi(M\cdot)=P\hatPhi$ on the one hand,
and the kernel $\ker L$ of the map
$$L : \C^r\times \zn \to  \C^r\times \zn
: (w_\alpha)_{\alp\in\zn} \mapsto (M^{|\alpha|}w_\alpha-
\sum_{0\leq \beta \leq \alpha}  
D^{\alpha - \beta}P(0)w_\beta:\  \alpha\in \zn),$$
on the other hand.
\end{thm}

\sl Proof. \rm Analogous to that of Theorem~\ref{charthm}. 
\eop

If the refinement equation is inhomogeneous, viz., a given function
is added to its right-hand side, then any solution to it is a sum of its 
specific solution and a solution to the corresponding homogeneous
refinement equation. This allows for generalizations of Theorems~\ref{charthm}
and~\ref{char2thm} to inhomogeneous equations as well. For the exact 
statement, see~\cite{JJS}.

\smallskip If $P(0)$ is ``regular'' in the sense that
$1$ is the largest dyadic eigenvalue of it,
the characterization of Theorem~\ref{charthm}  is much simpler: every right
$1$-eigenvector of $P(0)$ gives rise to a solution $\Phi\in R(P)$. However,
it is easy to generate examples when the largest dyadic eigenvalue of $P(0)$ is
$>1$: for example one can replace $P$ by $2P$. In that case, the solutions
in $R(2P)$ are obtained by differentiating suitably the solutions in $R(P)$.
We close this section with two results related to the current discussion.
In the first, we describe a general set up in which the solution space
$R(P)$ is decomposed into the sum of derivatives of solutions to ``regular''
refinement equations. In the second result, we provide an example when such
decomposition does not hold. Since the discussion here is somewhat
tangential to our main study of approximation orders, we skip the proof
of the following theorem.

\begin{thm} \label{layersthm} Given an $r\times r$ refinement
mask $P$,
let $R(P)$ denote the space of compactly supported
solutions to~(\ref{*}), and let $N$ be the maximal integer $n$ for
which $2^n\in \spec(P(0))$.
Suppose, that that we can find two $r\times r$ function-valued matrices
$T$ and $\tilP$ such that: (i) $T$ is analytic and
invertible around the origin, (ii) the entries of $\tilP$ are
trigonometric polynomials, (iii) the matrix $T(2\cdot)P-\tilP T$ has a
zero of order $N+1$ at the origin, and  (iv) the Taylor expansion of
degree
$N$ of $\tilP$ around the origin is block-diagonal and the spectrum of
each
block evaluated at zero intersects the set $\{2^{j}:\ j=0,\ldots,N\}$ at
no
more than one point.  Let $\Phi$ be in $R(P)$, and assume that each of the
entries of $\hatPhi$ has a zero of order $l$ at the origin. Then $\Phi$
admits a representation
\begin{equation}
\Phi=\sum_{j=l}^N p_j(D)\Phi_j,  \label{layers}  \end{equation}
with $\Phi_j\in R(P/2^j)$, $\hatPhi_j(0)\not=0$,  and  $p_j$ is a
homogeneous
polynomial of degree $j$, $j=l,\ldots,N$.

\end{thm}

As mentioned before, the solution space $R(P)$ does not always have such 
structure, as the following counterexample demonstrates.

\begin{exmpl} For some masks $P$, the decomposition~(\ref{layers})
from Theorem~\ref{layersthm} is not valid.
\end{exmpl}

\sl Proof. \/ \rm Let $d=2$ and let the mask $P$ satisfy the following
conditions:
$$P(0)=\left[\begin{array}{ccc} 1 & 0 & 0 \\ 0 & 2 & 0 \\ 0 & -1 & 4
\end{array} \right],
\quad (D^{(0,1)}P)(0)= \left[\begin{array}{ccc} 0 & 0 & 0 \\ 0 & 0 & 0 \\
0 & 0 & 1 \end{array} \right], \quad  (D^{(1,0)}P)(0)=0.$$
Let us show that the following inclusion fails:
\begin{equation}
 \spa[(D^{\alpha}\widehat{\Phi}(0))_{|\alpha|\leq N} :
 \Phi\in R(P), \widehat{\Phi}(0)=0]\subseteq \sum_{j=1}^d
\spa[( D^{\alpha}\left( ()^{e_j}\widehat{\Phi}\right)(0)  )_{|\alpha|\leq
N} :
\Phi\in R\left( {1\over 2} P \right)]. \label{intmed}
\end{equation}
Here $e_j$ is the vector in $\Z_+^d$ with $1$ in position $j$ and
zeros elsewhere. By Theorem~\ref{charthm}, this is equivalent to
the fact that~(\ref{layers}) fails.

The sequences $w$ in $\C^{r\times {N+ d \choose N-1}}$ (in our case
$r=3$, $d=N=2$) indexed by $\alpha$,  $|\alpha|<N$, we envision as `long'
vectors with the components $w_\alpha$, each of length $r$, all stacked
together in some fixed order, e.g., in the graded lexicographic order
of the $\alpha$'s.

The relation~(\ref{intmed}) is, again by Theorem~\ref{charthm}, equivalent
to the following:
\begin{equation}
 \ker L_0\subseteq \sum_{j=1}^d \ker L_j, \label{claim}
\end{equation}
where
$$L_0: \C^{r\times {N+ d \choose N-1}} \to  \C^{r\times {N+ d \choose
N-1}}
: (w_\alpha) \mapsto (2^{|\alpha|}I -P(0))w_\alpha-
\sum_{0< \beta < \alpha}
D^{\alpha - \beta}P(0)w_\beta, \qquad 1 \leq |\alpha|\leq N,$$
$L_j: \C^{r\times {N+ d \choose N-1}} \to  \C^{r\times {N+ d \choose N-1}}
$:
$$(w_\alpha) \mapsto \cases{ (2^{|\alpha|}I -P(0))w_\alpha-
\sum_{e_j\leq \beta < \alpha}
D^{\alpha - \beta}P(0)w_\beta, & if $\alpha \ge e_j$ \cr
w_\alpha & otherwise,  } \qquad 1 \leq |\alpha|\leq N.$$
Now,~(\ref{claim}) fails iff its dual statement
\begin{equation} \ran L_0^* \supseteq \cap_{j=1}^d
\ran  L_j^* \label{dualclaim} \end{equation}
 fails.  Here~(\ref{dualclaim}) is obtained from~(\ref{claim})
by taking orthogonal complements on both sides and using the property
$\ker A=(\ran A^*)^\perp$, which holds, in particular, for any linear
map acting on a finite-dimensional Hilbert space.

Now let $w_{(0,1)}\eqbd \left[ \begin{array}{c} 0 \\ 0 \\ 2 \end{array}
\right],
\qquad w_{(1,1)}\eqbd \left[ \begin{array}{c} 0 \\ 1 \\ -2 \end{array}
\right].$
Then $$w\eqbd  \left[ \begin{array}{c} w_{(0,1)} \\ 0 \\ 0 \\  0 \\ 0
\end{array} \right]=
L^*_{(0,1)}\left[ \begin{array}{c} 0 \\ 0 \\ 0 \\  w_{(1,1)} \\ 0
\end{array} \right]=
L^*_{(1,0)}\left[ \begin{array}{c} w_{(1,0)} \\ 0 \\ 0 \\  0 \\ 0
\end{array} \right]$$
but, by direct calculation, the vector  $w$ is not in the range of
$L^*_0$.
\eop

\subsection{Coherent approximation orders \label{cohntsec}}

The general theory of approximation orders of FSI spaces (\S4) focuses
on the individual space $S_\Phi$ and its properties. In contrast, when 
studying the solutions of the refinement equation~(\ref{*}), we believe 
that the focus should be on the interplay among the various solutions, 
in other words on their ``common ground''. An attempt to establish a theory 
that treats simultaneously all the solutions of (\ref{*}) should be done
with care: it is easy to show that different solutions of the same refinement
equation may have completely different properties, as the following
discussion makes clear.

\begin{dis}\label{disc} For $j=1,\ldots,r$, let $\phi_j$ be a (scalar-valued)
refinable function with (scalar) mask $p_j$. That is, each $p_j$ is
$2\pi$-periodic and $\hatphi_j(2\cdot)=p_j\hatphi_j$. Define
$P\eqbd \diag(p_1,\ldots,p_r)$, $\Phi\eqbd (\phi_j)_{j=1}^r$. Then $\Phi$ is
a  refinable vector with mask $P$. For each fixed $j$, the vector
$\Phi_j$ whose $j$th entry is $\phi_j$ and all other entries
are $0$ is refinable with respect to $P$. Since we may select the original
refinable elements $(\phi_j)$ in a completely arbitrary manner,
it is clear that the different solutions $(\Phi_j)_j$ to the same
refinement equation may be very different one from the other.
\end{dis}

This discussion reveals another difficulty that arises when dealing
with different solutions to the same refinement equation: with
$G_j$ the Gramian of $\Phi_j$, that Gramian is {\it{singular}} at the
origin. It is well-known that this is not an accident:

\begin{res}[{\cite{JS}}] Let $\Phi \subset \l2$ be a compactly
supported refinable vector with Gramian $G_\Phi$. If $G_\Phi(0)$ is 
invertible, then the spectral radius $\varrho(P(0))$ of $P(0)$ is 
equal to $1$, $1$ is the only eigenvalue on the unit circle, and
$1$ is a simple eigenvalue. \end{res}

That is, the Gramian of a refinable function is invertible at zero only 
if the spectrum of $P(0)$ is of a special nature, which, in particular, 
implies that the refinement equation has a {\it{unique}} solution.  
We note that the analysis of the approximation order of this case 
(viz., a refinable vector whose Gramian is invertible at the origin) 
is carried out in \cite{BDR4} and \cite{J1} and is {\it{not\/}}
among our objectives here (although we will recall those results momentarily).

In order to deal with all the solutions of a fixed refinement
equation in a combined fashion, we introduce first the notions of
the {\it{combined Gramian}} and the {\it{coherent approximation order}}
of the solutions. Let $P$ be a refinement mask, and let
$(\Phi_1,\ldots,\Phi_n)$ be a basis for the solution space $R(P)$ of the
underlying refinement equation (\ref{*}). Assuming that, for some
$s\in\R$ and for every $j=1,\dots,n$, $\Phi_j\subset \wsr$, we define
the {\bfi{combined Gramian}} $G_{R(P),s}$\index{$G_{R(P),s}$ -- combined
Gramian in $\wsr$} of the refinement equation
(\ref{*}) to be the {\it{sum}} of the individual Gramians:
$$G_{R(P),s}\eqbd \sum_{j=1}^n G_{\Phi_j,s}.$$
Although the above definition depends on the particular basis that
we choose for the solution space, our subsequent analysis of $G_{R(P),s}$
is independent of the basis' choice for the following reason.
Let $B=(\Phi_1,\ldots,\Phi_n)$ be a basis for $R(P)$. We consider
$B$ as an $r\times n$ matrix. Thanks to the identity
\begin{equation} \sum_{l=1}^n\widehat{\Phi_l}\widehat{\Phi_l}^*
=\hatB\hatB^\ast,\label{qform}
\end{equation}
we conclude that 
$$G_{R(P),s}=\sum_{\alp\in 2\pi\Zd}(\hatB \hatB^\ast
|\cdot|^{2s})(\cdot+\alp).$$
A new basis for $R(P)$ can be written as $BM$, with $M$ an $n\times n$ constant matrix.
Thus the combined Gramian for the new basis has the form
$$\tilG_{R(P),s}\eqbd \sum_{\alp\in 2\pi\Zd}(\hatB MM^\ast \hatB^\ast
|\cdot|^{2s})(\cdot+\alp).$$
Therefore, for some constants $c$, $C>0$, $$ c v^*\tilG_{R(P),s} v \leq 
v^* G_{R(P),s}v \leq C v^* \tilG_{R(P),s}v, $$
for any vector $v$. Using the above inequalities, one can easily  
check that all our subsequent results are independent of the choice of $B$.
We also use the notion of the truncated combined Gramian:\index{$G^0_{R(P),s}$ -- 
truncated combined Gramian in $\wsr$} 
$$G^0_{R(P),s}\eqbd \sum_{j=1}^n G^0_{\Phi_j,s}.$$

{\bf Definition\/.} \rm Let $P$ be a refinement mask whose solution
space $R(P)$ lies in $\wsr$. Let $G_{R(P),s}$ be the corresponding combined
Gramian. We say that $R(P)$ (or, in short, $P$) provides
{\bfi{coherent approximation order}} $k$ if the following condition holds:
there exists a neighborhood $\Omega$ of $0$ such that 
$$  \hbox{the function }
\Mu_{P,s,k} : \omega \mapsto {1\over |\omega|^{2k-2s}} \inf_v 
{v^* G^0_{R(P),s}(\omega) v \over v^* G_{R(P),s}(\omega) v } \quad
\hbox{\rm belongs to} \quad L_\infty(\Omega).$$  
Here, $G_{R(P),s}^0$ is the truncated version of the combined Gramian, i.e.,
the expression obtained after subtracting from $G_{R(P),s}$ the terms
indexed by $\alp=0$ (cf.\   the analogous
definition of the truncated Gramian from Theorem~\ref{fsisob}).

While Theorem \ref{fsisob} provides ample motivation for the above 
definition (specifically, it shows that the coherent approximation
order coincides with the usual approximation order notion in case
the solution space of (\ref{*}) is one-dimensional), we note that
the coherent notion of approximation order does not translate
immediately into any clear statement on the approximation order of the
individual solutions. 

\begin{dis} Let us continue with the example in Discussion \ref{disc}.
We observe that in the case discussed there,
$G_{R(P),s}=\diag(G_{\phi_1,s},\ldots,G_{\phi_n,s})$, with $G_{\phi_j,s}$ the 
(scalar) Gramian of $\phi_j$, i.e., $[\hatphi_j,\hatphi_j]$ in the $L_2$-case.
It follows easily then that the coherent approximation order matches or 
exceeds the approximation order provided by $\phi_j$ (for any value of 
$j$).  
\end{dis}

In order to advance our discussion, we consider vectors $v$ that 
realize the coherent approximation order $k$. That is,
with $\Ome\subset\Rd$ some neighborhood of the origin, 
$$\Ome\ni\ome\mapsto v(\ome)\in\C^r$$
is measurable, and,  a.e.\ on $\Ome$,
\begin{equation}  {v^*(\ome) G^0_{R(P),s}(\omega) v(\ome) \over v^*(\ome)
G_{R(P),s}(\omega) v(\ome) } =O(|\ome|^{2k-2s}).\label{supervector}
\end{equation}
We call such $v$ a {\bfi{universal supervector}} (of order $k$).
A vector $v$ is a {\bfi{regular universal supervector}} if (\ref{supervector}) 
can be replaced by the conditions that, near the origin,
\begin{equation}   {v^* G_{R(P),s} v\over v^\ast v}\sim |\cdot|^{2s}
\quad \hbox{and}\quad 
{v^* G^0_{R(P),s} v \over v^* v}=O(|\cdot|^{2k}).
\label{regsupervector}
\end{equation}
A regular universal supervector is clearly a universal supervector.

\begin{dis} The existence of a universal supervector is implied
(almost automatically) by the definition of coherent approximation order.
The proof of this fact parallels the proof of the superfunction existence 
(Theorem~\ref{super}) and is therefore omitted. The regularity of a universal 
supervector $v$ may be implied by either of the following two stronger assumptions:

(1) The combined Gramian $G_{R(P),s}$ is invertible a.e.\ around the origin, 
and the norm of its inverse there satisfies
$$\|G^{-1}_{P,s}\|=O(|\cdot|^{-2s}).$$ 
Indeed, in that case $v^* G_{P.s}v\geq c v^* v |\cdot|^{2s}$ 
for some positive constant $c$, since $\|G^{-1}_{P,s}\|$ 
is proportional to the reciprocal of the smallest eigenvalue 
$\sigm_\mi (G_{R(P),s})$ of $G_{R(P),s}$ and  $v^* G_{R(P),s}v \geq \sigm_\mi (G_{R(P),s}) v^* v$.
On the other hand, if  $v$ is a universal supervector, then
$$  v^* G_{p,s} v\leq \const v^*   \sum_{\Phi\in B}\hatPhi \hatPhi^* |\cdot|^{2s} v
\leq \const  |\cdot|^{2s}, $$ 
where the last inequality follows from the fact that all $\Phi\in B$ are 
compactly supported, so their Fourier transforms are bounded around the origin.
Therefore, the first, hence all the conditions in (\ref{regsupervector}) are 
satisfied.

(2) For one of the solutions
$\Phi\in R(P)$, $|v^\ast\hatPhi|/|v|\ge c>0$, a.e.\ in some neighborhood of the
origin. Indeed, then 
$$ v^*   \sum_{\Phi\in B}\hatPhi \hatPhi^* |\cdot|^{2s} v \sim |\cdot|^{2s} v^* v$$
and the conditions~(\ref{regsupervector}) follow from the fact that $v$ is
a universal supervector.

\end{dis}

We now connect among the notions of coherent orders, approximation 
orders, and regular universal supervectors. It is worthwhile to note 
that  the following result does not invoke the fact that
$R(P)$ comprises the solutions to (\ref{*}). We do not even need the
fact that the individual vectors in $R(P)$ are refinable.

\begin{thm} \label{coherent} Assume that the refinement mask $P$ 
provides coherent approximation order $k$ in $\wsr$. Then
\begin{description}
\item{(a)} Let $S_P\subset\wsr$ be the shift-invariant space generated by
$R(P)$ (i.e., it is the smallest closed shift-invariant subspace of 
$\wsr$ that contains each entry of each vector in $R(P))$. Then $S_P$ is 
an FSI space and provides approximation order $k$.
\item{(b)} Let $v$ be a regular universal supervector of order $k$
that is bounded in a neighborhood of the origin.
Let $\Phi$ be a solution of the refinement equation. Then
\begin{description}
\item{(i)} $v^\ast G^0_{\Phi,s}v=O(|\cdot|^{2k})$ around
the origin. 
In particular, the function $\psi$ defined by 
$\hatpsi\eqbd v^\ast\hatPhi$ satisfies the Strang-Fix conditions of order 
$k$.
\item{(ii)} If, for some positive $c$,
$|v^\ast\hatPhi|\ge c$ a.e.\ in some neighborhood of the origin,
then $S_\Phi$ provides approximation order $k$. Moreover, with $\psi\in
S_\Phi$ defined by $\hatpsi\eqbd v^\ast\hatPhi$, the PSI space $S_\psi$ 
already provides that approximation order.
\end{description}
\end{description}
\end{thm}

\sl Proof. \rm 
(a) Let $B$ be a basis for $R(P)$. Then $S_P=S_F$, with $F$ any vector that
contains all the entries from all the vectors $b\in B$. Hence $S_P$ is
FSI.

Now, let $v:[-\pi,\pi]^d\to \C^r$ be a vector that realizes the coherent 
approximation order $k$, i.e., a.e.\ on $[-\pi,\pi]^d$,
\begin{equation}  
{v^* G^0_{R(P),s} v \over v^* G_{R(P),s} v } \le c |\cdot|^{2k-2s}. \label{tmp}
\end{equation}  
If follows that for a.e. $\ome\in [-\pi,\pi]^d$, there exists $\Phi\in B$
such that
\begin{equation}  
{v^*(\ome) G^0_{\Phi,s}(\ome) v(\ome) \over v^*(\ome) G_{\Phi,s}(\ome) 
v(\ome) } \le c |\ome|^{2k-2s}.\label{temp} \end{equation}  
This allows us to represent $[-\pi,\pi]^d$ as the disjoint union of 
sets $\Ome_\Phi$, $\Phi\in B$ such that (\ref{temp}) holds for every 
$\Phi\in B$ and a.e.\ $\ome\in\Ome_\Phi$. We need, furthermore,
to ensure that these sets are {\it measurable}.
We argue the measurability as follows. First, since $v$ is measurable, 
so are the functions from the left-hand side of~(\ref{temp}). 
Therefore,  the function 
$$ f_\mi\;\; : \;\;\ome \mapsto 
\min_{\Phi\in B}  {v^*(\ome) G^0_{\Phi,s}(\ome) v(\ome) \over v^*(\ome) 
G_{\Phi,s}(\ome) v(\ome) } $$ is also measurable. Thus, once we define
$(\Omega_\Phi)$ by
$$ \Omega_\Phi\eqbd  \{\ome \in [-\pi,\pi]^d : f_\mi(\omega)={v^*(\ome) 
G^0_{\Phi,s}(\ome) v(\ome) \over v^*(\ome) G_{\Phi,s}(\ome) v(\ome) } \}, 
\qquad \Phi\in B,$$ 
we obtain the requisite measurability.

Now, let $\tau_\Phi$, $\Phi\in B$, be the $2\pi$-periodic
extensions of the characteristic functions of $\Ome_\Phi$, $\Phi\in B$. 
Defining $\Phi_0$ via its Fourier transform as follows:
$$\hatPhi_0\eqbd \sum_{\Phi\in B}\tau_\Phi\hatPhi,$$
we conclude from Corollary \ref{sobcont} that each of the entries of $\Phi_0$
lie in $S_{P,s}$.  Consequently, $S_{\Phi_0,s}\subset S_{P,s}$.
On the other hand, the
definition of $\Phi_0$ implies that, a.e.\ on $[-\pi,\pi]^d$,
\begin{equation}  
{v^* G^0_{\Phi_0,s} v \over v^* G_{\Phi_0,s} v } \le c |\cdot|^{2k-2s}.
\end{equation}  
This, in view of Theorem \ref{fsisob}, shows that $S_{\Phi_0}$ provides approximation 
order $k$ (in $\wsr$), {\it{a fortiori\/}} its superspace $S_P$ provides that order.

(b): 
The regularity of the supervector $v$ implies that, around the origin,
$${v^\ast G_{R(P),s}^0v \over v^\ast v}=O(|\cdot|^{2k})$$
Since $v$ is assumed bounded, we conclude that
$$v^\ast G_{R(P),s}^0v=O(|\cdot|^{2k}),$$
and therefore
$$v^\ast G_{\Phi,s}^0v=O(|\cdot|^{2k})$$
for {\it{every}} $\Phi\in R(P)$. This proves the first part of (b)(i),
while the second part follows from the fact that each of the summands
$v^\ast |\cdot+\alp|^{2s}(\hatPhi\hatPhi^\ast)(\cdot+\alp)\, v$ (that 
together make up $v^\ast G_{\Phi,s}^0\, v$) is nonnegative, hence has
to vanish to order $2k$ as well.

As to (b)(ii), the analysis above shows that the function
$\psi$ defined by $\hatpsi\eqbd  v^* \hatPhi$  satisfies 
$$[\hatpsi,\hatpsi]_s^0=[\hatpsi,\hatpsi]_s-|\hatpsi|^2 |\cdot|^{2s}=
O(|\cdot|^{2k})$$
(near the origin). Since we further assume here that $|\hatpsi|\ge c>0$
around the origin, we also conclude that $[\hatpsi,\hatpsi]_s\ge
c^2 |\cdot|^{2s}$ there. Thus, (\ref{statPSI}) of Theorem~\ref{psisob} holds, 
and that theorem implies that $S_\psi(\wsr)$ provides approximation order 
$k$.
%
\eop

The first part of Theorem~\ref{coherent} leads to the following
conclusion:

\begin{corol} Let $P$ be a refinement mask and let $S_P\subset \wsr$ be 
the corresponding SI space. If $P$ provides a coherent approximation order
$k$, then there exists $\psi\in S_P$ for which the PSI space $S_\psi\subset
\wsr$ provides approximation order $k$.
\end{corol}

{\bf Remark.\/} Note that the combined Gramian $G_{R,s}$ can be defined 
for any {\em finite-dimensional\/} space $R$ of distributional solutions 
to the refinement equation~(\ref{*}). Likewise, the notion of (regular) 
universal supervectors makes sense with respect to any such space $R$ .
The requirement that $R$ be the space of all compactly supported solutions 
actually plays no role in the results of this section. The only condition 
used is that, for each $\Phi\in R$, its Fourier transform $\hatPhi$ be 
bounded around the origin. Therefore, all results of this section are 
applicable to this more general setup.


\subsection{\label{unisuper}Universal supervectors and sum rules}

\subsubsection{Known results: singleton solutions in $\ld$}

The characterizations to-date of the approximation power of refinable 
vectors are confined to the $\l2$ setup, and assume, at a minimum, that
the Gramian of the (necessarily unique) compactly supported solution
is invertible at zero (as well as at several additional points). These
characterizations allow one to deduce the approximation order
provided by the refinable vector directly from the mask.
Two relevant notions in this context are {\it{Condition $Z_k$}}
and the {\it{sum rules}}. We begin with the definition of the former.

{\bf Definition.\/} \rm Given $k>0$, we say that
the refinement mask $P$ {\bfi{satisfies Condition\/}} $Z_k$
if there exists a vector $v$ of trigonometric polynomials 
such that,
for each $l\in E$, the vector
$v^*(2\cdot)P-\delta_{l,0}v^\ast$ has a zero or order $k$ at $\pi l$,
while $v(0)\not=0$.
Here,
\begin{equation}
E\eqbd \{0,1\}^d,
\label{defE}
\end{equation}\index{$E$ -- vertices of the unit cube in $\R^d$}
is the set of vertices of the $d$-unit cube.

\begin{res}[{\cite{BDR4}}, {\cite{J1}}] \label{nice}  Let $P$ be a
refinement mask, and assume that $\dim R(P)=1$.
Let $\Phi$ be the unique solution of~(\ref{*}), and
assume that $\Phi\subset\ld$ is compactly supported,
and that $G_\Phi$ is invertible at the origin.  
Then, for $k\in\N$:
\begin{description}
\item{I.} 
If $P$ satisfies Condition $Z_k$ then $S_\Phi$ provides
approximation order $k$.
\item{II.} If $S_\Phi$ provides approximation order $k$ and if 
$G_\Phi$
is invertible at each point of $\pi E$, then $P$ satisfies Condition 
$Z_k$.  \end{description}   \end{res}

We note that the compact support assumption on $\Phi$ in the above-quoted
result can be weakened: the essential needed information is about the
behavior of the Gramian around $E$. We refer to
\cite{BDR4} and \cite{JP} for more details.

Condition $Z_k$ is written on the Fourier domain. It can be equivalently
expressed on the ``space'' domain. The equivalent space-based 
formulations of Condition $Z_k$ are colloquially known as the {\it{sum 
rules}}.  We provide, for completeness, the two frequently used 
versions of these sum rules.  The second is taken from \cite{BDR4} (see also 
\cite{HSS}), while the first is borrowed from \cite{J1}.

\begin{res} \label{sumrules} Let $\Phi\subset \l2$
be compactly supported with trigonometric refinement mask $P$. Let
$v$ be a vector of trigonometric polynomials. Then the following 
conditions are equivalent:
\begin{description}
\item{(a)} $P$ satisfies Condition $Z_k$ with respect to the current $v$.
\item{(b)} The pair $v,P$ satisfies the 1st version of sum rules:
with $(v_\gam)$ and $(P_{\gam})$ the Fourier coefficients of 
$v$ and $P$, respectively,
\begin{eqnarray*}
\sum_{\sigma\in \Z^d} \sum_{\gamma\in \Z^d} v^*_{\sigma-\gamma}
P_{l+2\sigma} q(l+2\gamma) = 
2^{-d}\sum_{\gam\in\Zd}v^\ast_{-\gam}q(\gam), \qquad l\in 
E, \quad q\in \Pi_{k-1}.  \end{eqnarray*} 
(Note that $v_\gam$ is a vector, $P_\gam$ is a matrix, and $q(\gam)$ is a
scalar.)
\item{(c)} The pair $v,P$ satisfies the 2nd version of sum rules:
With $v^\alp=D^\alp v(0)$, $\alp\in\Zd_+$, the Taylor coefficients of
$v$ at the origin, we have:
$$\sum_{\bet\le\alp}2^{|\alp-\bet|}(v^{\alp-\bet})^\ast(D^\bet P)(\pi
l)=\delta_{l,0}(v^\alp)^\ast,\quad l\in E,\ |\alp|<k.$$

\end{description}
\end{res}

\let\sig\sigma
\let\del\delta

\sl Proof. \rm  The second version of the sum rules is equivalent to Condition 
$Z_k$, as seen by applying $D^\alp$ to $v^\ast(2\cdot)P-\delta_{l,0}v^\ast$,
expanding the first term by Leibniz' rule, and evaluating the result
at $\pi l$, $l\in E$. 

The equivalence of the
first version to Condition $Z_k$ can be argued as follows: first,
recall that for a finitely supported
$s\in \openC^{\Zd}$, its Fourier series $\hat s$
has a zero of order $k$ at the origin
iff $\Pi_{k-1}$ lies in the kernel of the functional 
$$\lam_s: q\mapsto (s\ast q)(0),$$ where
$s\ast q$ is either the semi-discrete convolution or the discrete convolution 
(the statement is true for each of the two choices) of the sequence $s$ 
and the polynomial $q$.
Let $m\in E$, and let $s_{m}$ be the (vector-valued)
Fourier coefficients of the function
$v^\ast(2\cdot)P(\cdot+\pi m)-\delta_{m,0}v^\ast.$
Thus, Condition $Z_k$ tells that, for every $q\in \Pi_{k-1}$ and for
every $m\in E$, 
$$0=\lam_{s_{m}}(q)=\sum_{\gam,\eta\in\Zd}
v^*_{-\gam} P_\eta q(2\gam+\eta)e^{\pi i
\eta\cdot m}-\delta_{m,0}\sum_{\gam\in\Zd}v^*_{-\gam} q(\gam).$$
Fixing some $l\in E$,
we can write $\eta=l'+l+2\sig$, for suitable $\sig\in\Zd$, 
and $l'\in E$. Thus,
$$0=\sum_{l'\in E}\sum_{\gam,\sig} v^*_{-\gam} P_{l'+l+2\sig}q(2\gam+2\sig+l'+l)
e^{\pi i(l'+l)\cdot m} -\delta_{m,0}\sum_{\gam}v^*_{-\gam} q(\gam).$$
Multiplying the two sides of the last display by $e^{-\pi i l\cdot
m}$ and summing over $m$, we obtain
$$0=\sum_{l'\in E}\sum_{\gam,\sig} v^*_{-\gam} P_{l'+l+2\sig}q(2\gam+2\sig+l'+l)
\sum_{m\in E} e^{\pi il'\cdot m} -\sum_{\gam}v^*_{-\gam} q(\gam).$$
Thus,
$$0=2^d\sum_{\gam,\sig} v^*_{-\gam} P_{l+2\sig}q(2\gam+2\sig+l)
-\sum_\gam v^*_{-\gam} q(\gam).$$
Replacing $\gam$ by $\gam-\sig$ finishes the proof.
\eop

\subsubsection{New results: multiple solutions in $\ws$}

Our analysis of the multiple solution case is based on drawing a
connection between universal regular supervectors on the one hand, 
and Condition $Z_k$ (together with its associated sum rules) on the 
other hand. 
This approach requires some limited regularity formulated in terms of 
the combined Gramian $G_{R,s}$ of the space $R$ of refinable distributions 
(see Theorem~\ref{supercondzk} for the precise assumption).

So, let $P$ be an $r\times r$ refinement mask and let $R$ be a 
finite-dimensional space of solutions to~(\ref{*}) lying in some Sobolev 
space: $R\subset \wsr$ (in the sense that every entry of every
vector in $R$ lies in $\wsr$). Let $k$ be a positive number. We consider 
vector-valued functions $v$, that together with $P$ and $R$ satisfy the
following assumptions:

\begin{assum}\label{532assump}
\begin{enumerate}
\item There exists a neighborhood $\Ome$ of the origin, such that, a.e.\ on
$\Ome$,
\begin{equation}
\|G_{R,s}(\cdot+\pi l)\|=\cases{O(1),&$l\in E\bks0$,\cr
                               O(|\cdot|^{2s}),&$l=0$.\cr}
\label{gbound} \end{equation}
\item The entries of $v$ are $2\pi$-periodic and measurable. Moreover, 
$v^\ast v$ is bounded, and bounded away from zero, around the origin.
\item For some $k>0$, and for every $\Phi\in R$, the function $\phi$ 
defined by $\hatphi\eqbd v^\ast\hatPhi$ satisfies condition (\ref{techn2}).
\end{enumerate}
\end{assum}
%

Note that Assumptions~\ref{532assump} are valid (regardless of the value 
of $k$) whenever $R$ contains only compactly supported solutions, and $v$ 
is a vector-valued trigonometric polynomial, provided that $v(0)\not=0$.

\begin{thm} 
\label{corolours} Let $P$ be an $r\times r$ refinement mask, let 
$R\subset\wsr$ be a finite-dimensional space of solutions to~(\ref{*}), and 
let $v$ be a vector-valued function, so that $P$, $R$ and $v$ satisfy 
Assumptions \ref{532assump} with respect to some $k>0$. If $P$ and $v$ 
satisfy Condition $Z_k$, then $v^\ast G_{R,s}^0 v$ has a zero of order $2k$ 
at the origin, and, in addition, for every $\Phi\in R$, 
$v^\ast\hatPhi-(v^\ast\hatPhi)(0)=O(|\cdot|^k)$, provided that 
$v^\ast\hatPhi$ is smooth at the origin.  In particular:
\begin{description}
\item{(a)} If $v^\ast G_{R,s}v\sim |\cdot|^{2s}$ a.e.\ around the origin, 
then $v$ is a universal supervector (with respect to $R$) of order $k$, and 
hence $R$ provides coherent approximation order $k$. This, for example, is 
the case if $G_{R,s}$ is invertible (a.e.) in a neighborhood of the origin, 
and $\|G^{-1}_{R,s}\|=O(|\cdot|^{-2s})$ a.e.\ there.
\item{(b)} If, for some $\Phi\in R$, $|v^\ast\hatPhi|\ge c>0$ a.e. 
around the origin, then (a) applies, and we further conclude that 
$S_\Phi$ provides approximation order $k$.
\end{description} \end{thm}

\sl Proof. \rm  
Let $\alp\in 2\pi\Zd\bks0$, and let $m\eqbd m(\alp)\ge 1$ be the smallest 
integer for which $\alp/2^m\not\in 2\pi\Zd$. Let $\Phi\in R$. We prove, by 
induction on $m(\alp)$, that $v^\ast\hatPhi(\cdot+\alp)=O(|\cdot|^k)$. For 
$m=1$, we choose $l\in E\bks0$ such that $2\pi l-\alp\in 4\pi\Zd$. Since
$$v^\ast\hatPhi(\cdot+\alp)=v^\ast P \left( {\cdot\over 2}+\pi l\right)
\hatPhi\left({\cdot+\alp\over 2}\right),$$
the claim follows from
the fact that $v^\ast P({\cdot\over 2}+\pi l)=O(|\cdot|^k)$. For $m>1$, we 
write
$$v^\ast\hatPhi(\cdot+\alp)=v^\ast\left({\cdot\over 2}\right) 
\hatPhi\left({\cdot+\alp\over 2}\right)
+(v^\ast P\left({\cdot\over 2}\right)-v^\ast\left({\cdot\over 2}\right))
\hatPhi\left({\cdot+\alp\over 2}\right).$$
By Condition $Z_k$, $v^\ast P(\cdot/2)-v^\ast(\cdot/2)=O(|\cdot|^k)$.
In addition, since $m(\alp/2)=m-1$, the induction hypothesis yields that
$v^\ast(\cdot/2) \hatPhi((\cdot+\alp)/2)=O(|\cdot|^k)$, too.

Now, fix $\Phi\in R$ and define $\hatpsi\eqbd v^\ast\hatPhi$. Since
$v$ is bounded and $\psi\in\ws$, we can invoke Corollary~\ref{sobcont} 
and conclude that $\psi\in S_\Phi(\ws)$.  Since $\psi$ satisfies (\ref{techn2}) (as 
stipulated in Assumptions \ref{532assump}), and since, by our argument above,
$\hatpsi(\cdot+\alp)=O(|\cdot|^k)$, for every $\alp\in 2\pi\Zd\bks0$,
we see that
\begin{equation} 
\sum_{\alp\in 2\pi\Zd\bks0}|\hatpsi(\cdot+\alp)|^2 |\cdot+\alp|^{2s}
=O(|\cdot|^{2k}).\label{temp1}
\end{equation}
However,  the left-hand-side in the above equality is
$v^\ast G^0_{\Phi,s}v$, and, hence, by summing
(\ref{temp1}) over a basis of $R$ we obtain that $v^\ast 
G^0_{R(P),s}v=O(|\cdot|^{2k})$.

Next, with $\psi$ as above, the case $l=0$ in Condition $Z_k$ together with
the boundedness of $\hatPhi$ around the origin (the latter is embedded
in Assumptions \ref{532assump}) imply that
$$\hatpsi(2\cdot)-\hatpsi=(v^\ast(2\cdot)P-v^\ast)\hatPhi=O(|\cdot|^k).$$
Once we assume $\hatpsi$ to be smooth around the origin, the above implies that
$\hatpsi-\hatpsi(0)=O(|\cdot|^k)$, as claimed.

The proofs of (a) and (b) are straightforward, hence are omitted.
\eop

\smallskip
Since the case of compactly supported solutions and a trigonometric polynomial 
$v$ is of central importance here, we record separately the statement of 
Theorem~\ref{corolours} for this case.

\begin{corol} \label{thmours} Let $P$ be a trigonometric polynomial refinement 
mask, let $R(P)\subset \ws$ be the space of all compactly supported solutions 
to the refinement equation~(\ref{*}), and let $v$ be a vector-valued 
trigonometric polynomial. Suppose that $P$ and $v$ satisfy Condition $Z_k$ 
for some $k>0$. Then:
\begin{enumerate}
\item For each $\Phi\in R(P)$, the (compactly supported) $\psi$ defined by
$\hatpsi\eqbd v^\ast\hatPhi$ satisfies the Strang-Fix conditions of order $k$, and in
addition, $\hatpsi-\hatpsi(0)=O(|\cdot|^k)$. Consequently, if 
$v^\ast(0)\hatPhi(0)\not=0$,
then $S_\Phi(\ws)$ provides approximation order $k$, and $\psi\in S_\Phi(\ws)$
is a corresponding superfunction.
\item If $v^\ast(0)\hatPhi(0)\not=0$ for {\sl some} $\Phi\in R(P)$,
then $P$ provides coherent approximation order $k$, and $v$ is a corresponding
universal regular supervector.
\end{enumerate}
\end{corol}

At present, we do not know whether Condition $Z_k$ is necessary for the 
provision of coherent approximation order $k$ in case $\dim R(P)>1$. The
results of this type that we are able to prove make strong assumptions
on the mask $P$. Below is one such result. The stringent assumption here
is that $P(0)=I$.  In what follows, we use the notation $G_{R(P),s}^0$ 
introduced in \S\ref{cohntsec} for the truncated combined Gramian.
We prove the result only for $s=0$, although it extends to other values 
of $s$ (at a cost of a few technical details and more awkward notation).

\begin{thm} \label{supercondzk} Let $P$ be an $r\times r$ 
trigonometric polynomial refinement mask and let $R(P)$ be the 
space of compactly supported solutions to~(\ref{*}). Suppose that 
$R(P)\subset \l2$ and that the combined Gramian $G_{R(P)}$ satisfies 
Assumption~\ref{532assump}.1, is smooth around each $l\in E$ and is 
boundly invertible around each $l\in E$. 
%
If $P(0)=I$, the following conditions are equivalent:
\begin{description}
\item{(a)} $P$ satisfies Condition $Z_k$ with some vector $v$
satisfying Assumptions~\ref{532assump}.2 and~\ref{532assump}.3.
\item{(b)} There exists a regular universal supervector $v$ of order $k$ 
for the space $R(P)$.
\end{description}
In addition, a regular universal supervector $v$ of order $k$ can
be always chosen so that, for every $\Phi\in R(P)$, 
$$v^\ast\hatPhi-(v^\ast\hatPhi)(0)=O(|\cdot|^k).$$
\end{thm}

\sl Proof. \rm In view of Theorem~\ref{corolours}, we only need to prove 
the implication (b$)\implies($a). 

We start the proof by noting the identities
\begin{equation}
G_{R(P)}(2\cdot)=
\sum_{l\in E}(PG_{R(P)}P^\ast)(\cdot+\pi l),
\label{superiden}
\end{equation}
and
\begin{equation}
G_{R(P)}^0(2\cdot)=
PG_{R(P)}^0P^\ast+\sum_{l\in E\bks0}(PG_{R(P)}P^\ast)(\cdot+\pi l).
\label{superidena}
\end{equation}
The first identity is straightforward (and is quite well-known;
cf.\ \cite{JS}).  The second one is obtained by the subtraction of
the identity
$$\sum_{\Phi\in B}(\hatPhi\hatPhi^\ast)(2\cdot)=
\sum_{\Phi\in B} P\hatPhi\hatPhi^\ast P^\ast$$
from the first one. Here, $B$ is the basis for $R(P)$ that was used to
define $G_{R(P)}$.

Let $\tau$ be any $2\pi$-periodic vector-valued function that satisfies the condition
\begin{equation}
\tau^\ast G_{R(P)}^0\tau=O(|\cdot|^{2k}), 
\label{tau}
\end{equation}
near the origin.  Then 
$(\tau^\ast G_{R(P)}^0\tau)(2\cdot)=O(|\cdot|^{2k})$ 
near the origin. Thus the evaluation at $\tau(2\cdot)$ of the quadratic form 
in the right-hand-side  of~(\ref{superidena}) leads to a function which has 
a zero of order $2k$ at the origin.  Since each summand there is 
nonnegative, it follows that, for every $l\in E\bks0$,
$$\tau^\ast(2\cdot)(PG_{R(P)}P^\ast)(\cdot+\pi l)\tau(2\cdot)=O(|\cdot|^{2k}).$$
However, $G_{R(P)}$ is assumed to be boundly invertible around
$\pi l$, hence we must have that $\tau^\ast(2\cdot)P(\cdot+\pi l)=O(|\cdot|^k)$,
near the origin, for every $l\in E\bks0$. 
In addition,
\begin{equation}
\tau^\ast(2\cdot)(PG_{R(P)}^0P^\ast)\tau(2\cdot)=O(|\cdot|^{2k}).
\label{taua}
\end{equation}

Now, let $v$ be a universal supervector. Then, (\ref{tau})
is satisfied for $\tau:=v$, hence $v$ satisfies 
the requirements in Condition $Z_k$ with respect to each $l\not=0$.
It remains to modify $v$ (if need be) so that
Condition $Z_k$ be satisfied at $l=0$, too. Note that so far
we have not used out special assumption on  $P$. Still, we already
know that (\ref{taua}) is satisfied for $\tau:=v$. 

In order to complete our argument, we assume that $P=I+O(|\cdot|^k)$ near 
the origin. We will revisit this condition after completing the main part 
of the proof. This additional condition, when applied to~(\ref{taua}) leads 
(once we take into account the boundedness and self-adjointness
of $G^0_{R(P)}$) to
\begin{equation}
\tau^\ast(2\cdot)G_{R(P)}^0\tau(2\cdot)=O(|\cdot|^{2k}).
\label{taub}
\end{equation}
Thus, we proved that (\ref{tau}) implies (\ref{taub}), and hence,
since (\ref{tau}) is satisfied for $\tau:=v$, we conclude that
\begin{equation}
v^\ast(2^n\cdot)G_{R(P)}^0v(2^n\cdot)=O(|\cdot|^{2k}) \qquad  \hbox{\rm for all } 
n\in \N.
\label{tauc}
\end{equation}
Our previous analysis then implies that
$v^\ast(2^n\cdot)P(\cdot+\pi l)=O(|\cdot|^k)$, for every
integer $n\ge 1$, and every $l\in E\bks0$.

Now, by forming a suitable finite linear combination of $v(2^n\cdot)$, 
$n=0,1,...$, we can construct a vector
$u$ such that $u(0)=v(0)$, while $u-u(0)=O(|\cdot|^k)$ at the origin.
Clearly, $u^\ast P(\cdot+\pi l)=O(|\cdot|^k)$, for every $l\in E\bks0$.
Since both $u-u(0)$ and $P-I$ have a $k$-fold zero at the origin,
we conclude that $u^\ast(2\cdot)P-u^\ast=O(|\cdot|^k)$. Thus
$u$ satisfies Condition $Z_k$.

We finally contend that there is no
loss of generality in the assumption that $P-I=O(|\cdot|^k)$ around 
the origin. Indeed, consider the transformation $P\mapsto T^{-1}(2\cdot) PT$,
where $T$ is a trigonometric polynomial-valued matrix,
such that $T(0)$ is invertible. For each $\ome$ in some small
neighborhood of the origin, we have the linear isomorphism defined
on  $\{\hatPhi(\omega):\ \Phi\in R(P)\}$ by 
$$\hatPhi(\omega)\mapsto  T(\omega)\hatPhi(\omega).$$
Since each $\hatPhi\in R(P)$ is entire, the isomorphisms induce
a corresponding one between the spaces $R(P)$ and $R( T^{-1}(2\cdot) PT)$.
Moreover, Condition $Z_k$ is invariant under this isomorphism,
since the vector $v^*(2\cdot)P-\delta_{l,0}v^*$ vanishes to order $k$
at $\pi l$ for each $l\in E$ if and only if so does the vector
$(v^*T)(2\cdot) (T^{-1}(2\cdot)PT) -\delta_{l,0} (v^*T)$.   

So, if we can show that we can choose a matrix-valued polynomial $T$ (of 
degree smaller than $k$) so that $T(0)$ is invertible and \begin{equation}
PT =T(2\cdot)+O(|\cdot|^k),\label{Pflat} \end{equation}
then our claim will follow. To this end, let $T(0)=I$ and let the 
derivatives $(D^\alpha T)(0)$ be defined inductively, according to
the partial order of multi-integers $\alpha$, as solutions
to the equation
\begin{equation}
\sum_{0\leq \beta\leq \alpha} (D^{\alpha-\beta} P)(0) (D^{\beta}T)(0)=
2^{|\alpha|} (D^\alpha T)(0),  \qquad 0<|\alpha|<k.  \label{Pexpand}
\end{equation}
This system is obtained by differentiating~(\ref{Pflat}) at the origin
and  is equivalent to~(\ref{Pflat}). For a fixed $\alpha$, the values
$(D^\beta T)(0)$, $\beta<\alpha$ are already chosen, 
and the coefficient of the term $D^\alpha T(0)$ is $2^{|\alpha|}-1\not=0$
(since $P(0)=I$).
Thus, (\ref{Pexpand}) 
has a solution $(D^\alpha T)(0)$. Thus, $P$ can be assumed to
be within $O(|\cdot|^k)$ of the origin. This completes the proof.
\eop

\subsubsection{Coherent polynomial reproduction}


We restrict our attention again to the space $R(P)$ of compactly
supported solutions to the refinement equation~(\ref{*}). We show 
that universal supervectors for $R(P)$ are also ultimately 
connected with polynomial reproduction using the shifts of any
compactly supported solution from $R(P)$. In short, we show that 
universal supervectors provide universal polynomial reproduction schemes:


\begin{thm} Let $P$ be a refinement mask whose space of compactly supported 
solutions $R(P)$ lies in $\wsr$. Let $v$ be a vector-valued trigonometric 
polynomial such that, for some $k>0$, {\rm any one\/} of the following 
conditions holds:
\begin{enumerate}
\item $v$ satisfies Condition $Z_k$.
\item $v$ and $P$ satisfy the 1st version of the sum rules
\item $v$ and $P$ satisfy the 2nd version of the sum rules
\item $v$ is a regular universal supervector of order $k$.
\end{enumerate}
Let $a$  be the vector-valued sequence of the Fourier coefficients of $v^*$,
and let $\Phi\in R(P)$. Then, with $a_1,\ldots,a_r$
the entries of $a$ and $\phi_1,\ldots,\phi_r$ the entries of
$\Phi$, the map
$$T_\Phi: q\mapsto \sum_{i=1}^r\phi_i\ast'(a_i \ast' q)=:\Phi\ast'(a\ast'q)$$
maps $\Pi_{<k}$ into itself.  The map is surjective (hence degree preserving)
if and only if $v^\ast(0)\hatPhi(0)\not=0$.
\end{thm}

\sl Proof. \rm By Result~\ref{sumrules}, conditions 1 through 3 are equivalent
one to the other,
and, by Theorem~\ref{corolours}, each of them implies 4. 
Condition 4, in turns, implies that,
for any $\Phi\in R(P)$, the compactly supported
function $\psi$ defined by $\hatpsi\eqbd v^* 
\hatPhi$ satisfies the Strang-Fix conditions of order $k$. The discussion 
preceding Theorem~\ref{repthm} now implies that the semi-discrete 
convolution operator $\psi\ast'$ reproduces polynomials of total degree at 
most $k-1$ and, moreover, preserves the degree if $\psi$ is nondegenerate, 
i.e., $\hatpsi(0)\neq 0$.  But $\psi$ itself is nothing but the semidiscrete 
convolution $\Phi\ast'\hat{v}^*$.  Since the convolution $\ast'$ is 
associative, we see that the function $\Phi\ast'(\hat{v}^* \ast' q)$ is a 
polynomial in $\Pi_{k-1}$ whenever $q\in \Pi_{k-1}$; moreover, it has 
exactly the same degree as $q$  whenever  $v^*(0)\hatPhi(0)=\hatpsi(0)\neq 0$.
\eop

\section*{Acknowledgments} We thank Carl de Boor for his critical 
reading of an earlier draft of this paper and Gautam Bharali for
a very helpful discussion of our earlier approach to section 4.9.

\printindex

\end{document}